\newif\ifjota

\jotafalse

\RequirePackage{fix-cm}
\ifjota
	\documentclass[smallextended,envcountsect]{svjour3}
\else
	\documentclass[smallextended,envcountsect]{svjour3_arxiv}
\fi
\smartqed

\usepackage{graphicx}
\usepackage{mathptmx}      % use Times fonts if available on your TeX system
\usepackage{amsmath,amssymb}
\usepackage{array}
\usepackage{enumitem}
\usepackage{xcolor}
\usepackage{bm}
\usepackage{multirow}
\usepackage{cite}

\usepackage{mathrsfs}

\providecommand{\ie}		{\emph{i.e\@.}\xspace}
\providecommand{\eg}		{\emph{e.g\@.}\xspace}

\providecommand{\myurl}[1][]	{\texttt{web.eecs.umich.edu/$\sim$fessler#1}\xspace}

\providecommand{\onweb}[1]	{Available from \myurl.}

\long\def\comment#1{}

\providecommand{\bcent}		{\begin{center}}
\providecommand{\ecent}		{\end{center}}
\providecommand{\benum}		{\begin{enumerate}}
\providecommand{\eenum}		{\end{enumerate}}
\providecommand{\bitem}		{\begin{itemize}}
\providecommand{\eitem}		{\end{itemize}}
\providecommand{\bvers}		{\begin{verse}}
\providecommand{\evers}		{\end{verse}}
\providecommand{\btab}		{\begin{tabbing}}	% use \= and \>
\providecommand{\etab}		{\end{tabbing}}
\newcounter{blist}
\providecommand{\blistmark}	{\makebox[0pt]{$\bullet$}}
\providecommand{\blistitemsep}	{0pt}
\providecommand{\blist}[1][]	{%
\begin{list}{\blistmark}{%
\usecounter{blist}%
\setlength{\itemsep}{\blistitemsep}%
\setlength{\parsep}{0pt}%
\setlength{\parskip}{0pt}%
\setlength{\partopsep}{0pt}%
\setlength{\topsep}{0pt}%
\setlength{\leftmargin}{1.2em}%
\setlength{\labelsep}{0.5\leftmargin}% centers label in margin
\setlength{\labelwidth}{0em}%
#1}% user over-rides
}
\providecommand{\elist}		{\end{list}}
\providecommand{\blistitemsep}	{0pt}
\providecommand{\bjfenum}[1][]	{%
\begin{list}{\bcolor{\arabic{blist}.} }{%
\usecounter{blist}%
\setlength{\itemsep}{\blistitemsep}%
\setlength{\parsep}{0pt}%
\setlength{\parskip}{0pt}%
\setlength{\partopsep}{0pt}%
\setlength{\topsep}{0pt}%
\setlength{\leftmargin}{0.0em}%
\setlength{\labelsep}{1.0\leftmargin}% label in margin?
\setlength{\labelwidth}{0pt}%
#1}% user over-rides
}

\newcounter{blistAlph}
\providecommand{\blistAlph}[1][]
{\begin{list}{\makebox[0pt][l]{\Alph{blistAlph}.}}{%
\usecounter{blistAlph}%
\setlength{\itemsep}{0pt}\setlength{\parsep}{0pt}%
\setlength{\parskip}{0pt}\setlength{\partopsep}{0pt}%
\setlength{\topsep}{0pt}%
\setlength{\leftmargin}{1.2em}%
\setlength{\labelsep}{1.0\leftmargin}% letter at left side of margin
\setlength{\labelwidth}{0.0\leftmargin}#1}%
}

\newcounter{blistRoman}
\providecommand{\blistRoman}[1][]
{\begin{list}{\Roman{blistRoman}.}{%
\usecounter{blistRoman}%
\setlength{\itemsep}{0.5em}\setlength{\parsep}{0pt}%
\setlength{\parskip}{0pt}\setlength{\partopsep}{0pt}%
\setlength{\topsep}{0pt}%
\setlength{\leftmargin}{4em}%
\setlength{\labelsep}{0.4\leftmargin}
\setlength{\labelwidth}{0.6\leftmargin}#1}%
}

\usepackage{bbm} % black-board math
\providecommand{\jfbbm}[1]	{\xmath{\mathbbm{#1}}} % needs bbm package
\providecommand{\qed}[1][0pt]	{\hfill\raisebox{#1}{\inmath{\Box}}} % need latexsym package
\providecommand{\reals}		{\jfbbm{R}}

\providecommand{\inprod}[2]	{\xmath{\mathop{\langle #1,\, #2 \rangle}\nolimits}}

\let\equivsave\equiv
\def\equiv{\xmath{\equivsave}}

\providecommand{\ba}[1]		{\left[ \begin{array}{#1}}
\providecommand{\ea}		{\end{array} \right]}
\providecommand{\be}		{\begin{equation}}
\providecommand{\ee}[1]		{\label{#1}\end{equation}}
\providecommand{\bea}		{\begin{eqnarray}}
\providecommand{\eea}[1]	{\label{#1}\end{eqnarray}}
\providecommand{\beas}		{\begin{eqnarray*}}
\providecommand{\eeas}		{\end{eqnarray*}}
\providecommand{\beals}[1][1]	{\begin{alignat*}{#1}}	% requires amsmath
\providecommand{\eeals}		{\end{alignat*}}

\providecommand{\berr}[2]{
\bgroup
\renewcommand{\theequation}{#1}
\be
#2
\ee{e,#1}
\egroup
\ignorespaces
}

\providecommand{\bearr}[2]{
\bgroup
\renewcommand{\theequation}{#1}
\bea
#2
\eea{e,#1}
\egroup
\ignorespaces
}

\providecommand{\inmath}	{\ensuremath}
\providecommand{\xmath}[1]	{\inmath{#1}\xspace}
\providecommand{\bmath}[1]	{\xmath{\bm{#1}}}	% needs \usepackage{bm}

\providecommand{\paren}[1]	{\xmath{\left(#1\right)}}

\providecommand{\abs}[1]	{\xmath{\left| #1 \right|}}

\providecommand{\xfrac}[2]	{\xmath{\frac{#1}{#2}}}

\providecommand{\Frac}[2]	{\xmath{{#1}/{#2}}}
\providecommand{\pfrac}[3][]	{\xmath{\!\paren{#1\frac{#2}{#3}}}}

\providecommand{\sinc}		{\Ofun{\mathrm{sinc}}}

\usepackage{jf-accent}
\usepackage{jf-fun}
\usepackage{jf-names}

\newcommand{\st} {\xmath{\text{s.t.}\:}}

\definecolor{mygreen}{rgb}{0.2,0.8,0.4}

\newtheorem{assumption}{Assumption}

\newcommand{\OGMG} {OGM\nobreakdash-G}

\newcommand{\symm} {\jfbbm{S}}

\newcommand{\cF} {\mathcal{F}_L(\reals^d)}
\newcommand{\cC} {C_L^{1,1}(\reals^d)}
\newcommand{\Xs} {\xmath{X_*(f)}}

\newcommand{\Bd} {\mathcal{B}}

\renewcommand{\P} {\operatorname{P}}

\newcommand{\zero} {\bmath{0}}

\newcommand{\A} {\bmath{A}}
\newcommand{\B} {\bmath{B}}
\newcommand{\C} {\bmath{C}}

\newcommand{\G} {\bmath{G}}
\newcommand{\Z} {\bmath{Z}}
\renewcommand{\S} {\bmath{S}}

\renewcommand{\aa} {\bmath{a}}
\newcommand{\bb} {\bmath{b}}

\renewcommand{\u} {\bmath{u}}

\newcommand{\h} {\bmath{h}}
\newcommand{\x} {\bmath{x}}
\newcommand{\y} {\bmath{y}}
\newcommand{\z} {\bmath{z}}
\newcommand{\g} {\bmath{g}}

\renewcommand{\del} {\bmath{\delta}}

\newcommand{\nnu} {\bmath{\nu}}

\newcommand{\tth} {\xmath{\tilde{h}}}
\newcommand{\ttheta} {\xmath{\tilde{\theta}}}

\newcommand{\hhh} {\xmath{\hat{h}}}
\newcommand{\htheta} {\xmath{\hat{\theta}}}

\ifjota
	\journalname{JOTA}
\else
	\usepackage{fullpage}
\fi

\usepackage{hyperref}
\hypersetup{
        linktoc=page,
        colorlinks=true,
        linkcolor=blue,
        citecolor=red,
        filecolor=magenta,
        urlcolor=cyan
}

\begin{document}

\title{Optimizing the Efficiency of First-Order Methods
for Decreasing the Gradient of Smooth Convex Functions
%\thanks{Communicated by Alexander Mitsos}
%\thanks{This work was supported in part by the National Research Foundation of Korea (NRF) grant
%funded by the Korea government (MSIT) (No. 2019R1A5A1028324),
%and the POSCO Science Fellowship of POSCO TJ Park Foundation.}
}

\titlerunning{Optimizing the Efficiency of First-Order Methods for Decreasing the Gradient}

\ifjota
	\author{
	Donghwan Kim
	\and
	Jeffrey A. Fessler
	}

	\institute{Donghwan Kim, Corresponding author \at
        %KAIST, 
	Korea Advanced Institute of Science and Technology (KAIST), \at
        Daejeon, Republic of Korea \at
        donghwankim@kaist.ac.kr
        \and
        Jeffrey A. Fessler \at
        University of Michigan, \at
        Ann Arbor, Michigan \at
        fessler@umich.edu
	}

	\date{Received: date / Accepted: date}

\else
	\author{
	Donghwan Kim$^1$
	\and
	Jeffrey A. Fessler$^2$
	}

	\institute{
	$^1$
        Department of Mathematical Sciences, KAIST, Republic of Korea \at
	$^2$
        Department of Electrical Engineering and Computer Science, 
	University of Michigan, USA \at
        \email{donghwankim@kaist.ac.kr, fessler@umich.edu}
	}

	\date{Date of current version: \today}
\fi

\maketitle

\begin{abstract}
This paper optimizes the step coefficients of first-order methods
for smooth convex minimization
in terms of the worst-case convergence bound (\ie, efficiency)
of the decrease in the gradient norm.
This work is based on the performance estimation problem approach. %~\cite{drori:14:pof}.
The worst-case gradient bound 
of the resulting method
is optimal up to a constant
for large-dimensional smooth convex minimization problems, %~\cite{nemirovsky:92:ibc};
under the initial bounded condition on the cost function value.
This paper then illustrates that
the proposed method %, named OGM-G, 
has a computationally efficient form
that is similar to the optimized gradient method. %(OGM). %~\cite{kim:16:ofo}.
\end{abstract}

\section{Introduction}

Large-dimensional optimization problems arise
in various modern applications
of signal processing, machine learning, control, communication, and many other areas.
First-order methods
are widely used for solving such large-scale problems
as their iterations involve only 
function/gradient calculations and simple vector operations.
However, they can require many iterations to achieve the given accuracy level.
Therefore, developing efficient first-order methods
has received great interest,
which is the main motivation of this paper.
In particular, this paper targets the decrease in the gradient 
for smooth convex minimization,
under the initial bounded condition on the cost function value.
This paper uses the performance estimation problem (PEP) in~\cite{drori:14:pof}
and constructs a new method called \OGMG. 

Among first-order methods for smooth convex minimization, 
Nesterov's fast gradient method (FGM)
\cite{nesterov:83:amf,nesterov:04}
has been used widely
because its worst-case \emph{cost function} inaccuracy bound (\ie, the cost function efficiency)
is optimal up to a constant,
under the initial bounded \emph{distance} condition
\cite{nesterov:04,nemirovsky:92:ibc}.
Recently, the optimized gradient method (OGM)
\cite{kim:16:ofo} 
(that was numerically first identified in~\cite{drori:14:pof} using PEP)
has been found to exactly achieve
the optimal worst-case rate
of decreasing the smooth convex cost functions~\cite{drori:17:tei},
leaving no room for improvement in the worst-case.
On the other hand, first-order methods
that decrease the \emph{gradient} at an optimal rate in~\cite{nemirovsky:92:ibc}
are yet unknown, even up to a constant.
The proposed \OGMG~method has such an optimal rate
under the initial bounded \emph{function} condition. % on the cost function value.
After the initial version of this paper was posted online~\cite{kim:18:otc},
a simple method using 
\OGMG~was constructed in~\cite{nesterov:20:pda}
that also has an optimal rate
under the initial bounded \emph{distance} condition.

Gradient rate analysis is useful 
both in theory %\footnote{
%The gradient rate analysis on smooth convex minimization in this paper
%can be possibly useful for the analysis on nonconvex optimization.
%}
(\eg, for a dual approach~\cite{nesterov:12:htm}
and a matrix scaling problem~\cite{allenzhu:18:htm-nips})
and in practice
(\eg, can be used as a stopping criterion).
In addition, unlike smooth convex minimization, a worst-case \emph{gradient} inaccuracy 
and an initial bounded \emph{function} condition are standard choices
for analyzing gradient methods for smooth \emph{nonconvex} minimization~\cite{drori:20:tco}.
Therefore, this work can provide a step towards better understanding
the convergence behavior of gradient methods for nonconvex minimization.

There is recent interest in developing accelerated methods
for decreasing the gradient (in convex minimization)
\cite{nesterov:12:htm,allenzhu:18:htm-nips,carmon:19:lbf,kim:18:ala,kim:18:gto}.
The best known worst-case gradient rate
is achieved by
FGM with a regularization technique in~\cite{nesterov:12:htm}
that is optimal up to a logarithmic factor.
However, a practical limitation of that method is 
that it requires knowledge of a bound on a value
such as the distance between the initial and optimal points.
In~\cite{kim:18:gto} we used PEP
to derive efficient first-order methods
that do not need knowledge of such unavailable values.
However, the methods in~\cite{kim:18:gto} are far from achieving the optimal rate
(not even up to a logarithmic factor),
due to strict relaxations introduced to PEP in~\cite{kim:18:gto}.
The methods in~\cite{nesterov:12:htm,kim:18:ala,kim:18:gto,ghadimi:16:agm,monteiro:13:aah}
also achieve a similar nonoptimal rate. %(but with larger constant).
Thus, 
there is still room to improve
the worst-case gradient convergence bound
of the first-order methods
for smooth convex minimization.

This paper optimizes the step coefficients of first-order methods
in terms of the worst-case gradient decrease
using PEP~\cite{drori:14:pof,taylor:17:ssc},
%leading to a new method called
yielding \OGMG.
The new analysis
avoids the (unnecessary) strict relaxations on PEP in~\cite{kim:18:gto}.
This paper then shows that
\OGMG~has an equivalent efficient form that is similar to OGM,
and thus has an inexpensive per-iteration computational complexity.
\OGMG~attains the optimal bound of the worst-case gradient norm
up to a constant
under the initial bounded \emph{function} condition~\cite{nemirovsky:92:ibc}.
On the way, this paper also provides a new exact worst-case gradient bound
for the gradient method (GM).

%{\cblue
The initial bounded condition on the distance between initial and optimal points is a standard assumption,
whereas the initial condition on the cost function value of interest in this paper 
is less popular.
However, sometimes a constant for the latter bounded condition is known,
while a constant for the former condition is either not known or difficult to compute,
making the latter condition more useful.
In addition, there are cases where the latter initial bounded function condition holds,
but the former condition does not.
One such example is an unregularized logistic regression of an overparameterized model for separable datasets
\cite{nacson:19:cog,soudry:18:tib},
which does not have any finite minimizer.
Therefore, this paper's analysis under the initial bounded condition 
has value for such cases.
%}

%
Section~\ref{sec:prob} reviews a smooth convex problem and first-order methods. 
Section~\ref{sec:eff} reviews 
the efficiency of first-order methods and its lower bound. 
%in~\cite{nemirovsky:92:ibc,drori:17:tei}.
Section~\ref{sec:pep} studies the PEP approach~\cite{drori:14:pof}
and provides relaxations for analyzing the worst-case gradient decrease.
Section~\ref{sec:gm,ic2}
uses the relaxed PEP to provide the exact worst-case gradient bound for GM.
Section~\ref{sec:optpep}
optimizes the step coefficients of the first-order methods
using the relaxed PEP,
and develops an efficient first-order method
named \OGMG~under the initial function condition.
Section~\ref{sec:conc} concludes the paper.

\section{Problems and Methods}
\label{sec:prob}

\subsection{\bf Smooth Convex Problems}

We are interested in efficiently 
solving the following smooth and convex minimization problem:
\begin{align}
f_* := \inf_{\x\in\reals^d} f(\x)
\tag{M}
\label{eq:prob}
,\end{align}
where we assume that 
%\begin{itemize}[leftmargin=2em]
%\item 
the function $f\;:\;\reals^d\to\reals$
is a convex function of the type $\cC$, \ie,
its gradient $\nabla f(\x)$ is Lipschitz continuous:
\begin{align}
||\nabla f(\x) - \nabla f(\y)|| \le L||\x - \y||,
\quad \forall \x,\y\in\reals^d
\label{eq:Lsmooth}
\end{align}
with a Lipschitz constant $L>0$, 
where $||\cdot||$ denotes the standard Euclidean norm.
%and
%\item the optimal set $\Xs := \argmin{\x\in\reals^d}f(\x)$ is nonempty.
%\end{itemize}
\begin{definition}
The class of smooth convex functions satisfying the two above conditions
is denoted by $\cF$. 
\end{definition}
\begin{definition}
The optimal set of $f$ is defined by
\begin{align}
\Xs := \argmin{\x\in\reals^d}f(\x) = \{\x\in\reals^d\;:\;f(\x) = f_*\}
.\end{align}
\end{definition}
%is nonempty.}
%
We further assume one of the following two initial conditions,
where the latter is especially useful
when~\eqref{eq:prob} does not have a finite minimizer,
\ie, $\Xs = \emptyset$.

\begin{assumption}[IFC]
%We assume that 
The set $\Xs$ is nonempty,
and %the difference between an initial function value and the optimal function value
	%at a point $\x_*\in\Xs$
	%is bounded as
        an initial point $\x_0$ satisfies
	\begin{align}
	f(\x_0) - f_* \le \frac{1}{2}LR^2 
	\quad \text{for a constant $R>0$.}
	\tag{IFC}
	\label{eq:init,func}
	\end{align}
\end{assumption}
%We also consider a condition that assumes that
\begin{assumption}[IFC$'$]
%The difference between an initial function value and the function value after $N$ iteration
%	of a given method
%        is bounded as
An initial point $\x_0$ and the $N$th iterate $\x_N$ of a given method satisfy
        \begin{align}
        f(\x_0) - f(\x_N) \le \frac{1}{2}LR_N^2
        \quad \text{for a constant $R_N>0$.}
        \tag{IFC$'$}
        \label{eq:init,func,N}
        \end{align}
\end{assumption}
Note that $f(\x_0) - f(\x_N) \le f(\x_0) - f_*$ for any $\x_N$.

\subsection{\bf First-Order Methods}

To solve a large-dimensional problem~\eqref{eq:prob},
we consider first-order methods
that iteratively gain first-order information,
\ie, values of the cost function $f$ and its gradient $\nabla f$
at any given point in $\reals^d$.
The computational effort for acquiring those values 
depends mildly on the problem dimension.
We are interested in developing a first-order method
that efficiently generates a point $\x_N$ after $N$ iterations
(starting from an initial point $\x_0$)
that minimizes the worst-case absolute gradient inaccuracy
under the initial function condition~\eqref{eq:init,func}.
\begin{definition}
The gradient efficiency is defined as the worst-case absolute gradient inaccuracy
\begin{align}
\sup_{f\in\cF} ||\nabla f(\x_N)||^2
\label{eq:epsgrad}
.\end{align}
\end{definition}

For simplicity in Sects.~\ref{sec:pep},~\ref{sec:gm,ic2} and~\ref{sec:optpep} 
that use the PEP approach
(as in~\cite{drori:14:pof}),
we consider the following \emph{fixed-step} first-order methods (FSFOM):
\begin{align}
\x_{i+1} = \x_i - \frac{1}{L}\sum_{k=0}^{i}h_{i+1,k}\nabla f(\x_k)
	\quad i=0,\dots,N-1,
\label{eq:fsfom}
\end{align}
where $\h := \{h_{i+1,k}\} \in \reals^{N(N+1)/2}$
is a tuple of fixed-step coefficients
that do not depend on %$f$ and $\x_0$ ({\cmag and} $R$ {\cmag or} $R_N$). 
$f$, $\x_0$ and $R$ (or $R_N$). 
This FSFOM class 
includes (fixed-step) GM 
%$\x_{i+1} = \x_i - \frac{1}{L}h_{i+1,i}\nabla f(\x_i)$, 
(\ie, $h_{i+1,k} = 0$ for $k<i$),
(fixed-step) FGM \cite{nesterov:83:amf,nesterov:04} (see~\cite{drori:14:pof}), 
OGM \cite{kim:16:ofo}, and the proposed \OGMG,
but excludes line-search approaches,
such as 
a backtracking version of FGM in~\cite{beck:09:afi} and
an exact line-search version of OGM in~\cite{drori:19:efo}.

\section{Efficiency of First-Order Methods}
\label{sec:eff}

This paper seeks to improve the \emph{efficiency} 
of first-order methods,
where the efficiency consists of the following two parts:
the computational effort for selecting a search point 
(\eg, computing $\x_{i+1}$ in~\eqref{eq:fsfom} given $\x_i$ and $\{\nabla f(\x_k)\}_{k=0}^i$),
and the number of evaluations of the cost function value and gradient
at each given search point
to reach a given accuracy.
This paper considers both parts of the efficiency,
particularly focusing on the latter part, 
as also detailed in this section.
Regarding the former aspect of the efficiency,
we later show that the proposed method has an efficient form,
similar to (fixed-step) FGM and OGM,
requiring computational effort comparable to that of a (fixed-step) GM.

An efficiency estimate of an optimization method is
defined by the worst-case absolute inaccuracy.
One popular choice of the worst-case absolute inaccuracy
is the worst-case absolute cost function inaccuracy.
\begin{definition}
The cost function efficiency is defined as the worst-case absolute cost function inaccuracy
\begin{align}
\sup_{f\in\cF} f(\x_N) - f_*
\label{eq:epsfunc}
.\end{align}
\end{definition}
When analyzing the cost function efficiency, we usually 
consider the following initial condition.
\begin{assumption}[IDC]
	The set $\Xs$ is nonempty,
	and %the distance between initial and optimal points
	%are bounded as
	an initial point $\x_0$ satisfies
	\begin{align}
	||\x_0 - \x_*|| \le \overline{R}
	\quad \text{for a constant $\overline{R}>0$,} 
	\tag{IDC} 
	\label{eq:init,dist}
	\end{align}
	for some $\x_* \in \Xs$.
\end{assumption} 
Under~\eqref{eq:init,dist},
GM has an $O(1/N)$ cost function efficiency~\eqref{eq:epsfunc}
\cite{nesterov:04},
and this rate was improved to $O(1/N^2)$ rate
by FGM~\cite{nesterov:83:amf,nesterov:04}.
This efficiency was further optimized by OGM~\cite{drori:14:pof,kim:16:ofo},
which was shown to exactly achieve the optimal efficiency in~\cite{drori:17:tei}.

Compared to the worst-case \emph{cost function} inaccuracy~\eqref{eq:epsfunc},
the worst-case absolute \emph{gradient} inaccuracy~\eqref{eq:epsgrad}
has received less attention
\cite{nemirovsky:92:ibc,nesterov:12:htm,taylor:17:ssc,taylor:18:ewc}.
%this paper optimizes this gradient efficiency~\eqref{eq:epsgrad}.
%
For the initial \emph{distance} condition~\eqref{eq:init,dist},
GM has an $O(1/N^2)$ gradient efficiency
\cite{nesterov:12:htm},
while FGM with a regularization technique~\cite{nesterov:12:htm}
that requires the knowledge of (practically unavailable) $\overline{R}$
achieves $O(1/N^4)$ up to a logarithmic factor.
This is the best known rate,
where the rate $O(1/N^4)$ is the optimal gradient efficiency for given~\eqref{eq:init,dist}
\cite{nemirovsky:92:ibc}.
On the other hand,
the papers~\cite{nesterov:12:htm,kim:18:ala,kim:18:gto,ghadimi:16:agm,monteiro:13:aah}
studied first-order methods that do not require knowing $\overline{R}$ 
and that have $O(1/N^3)$ gradient efficiency,
but none of them (including~\cite{nesterov:12:htm}) have the optimal efficiency
(even up to a constant).

On the other hand,
gradient efficiency with the initial \emph{function} condition~\eqref{eq:init,func}
has received even less attention
\cite{nemirovsky:92:ibc,taylor:18:ewc}.
It is known to have $O(1/N^2)$ optimal efficiency~\cite{nemirovsky:92:ibc}.
Section~\ref{sec:gm,ic2} provides the exact $O(1/N)$ rate of GM,
which was studied numerically
for the more general nonsmooth composite convex problems in~\cite{taylor:18:ewc}.
The paper~\cite{carmon:19:lbf} discusses 
that FGM with a regularization technique~\cite{nesterov:12:htm}
with~\eqref{eq:init,func}
also achieves the optimal worst-case gradient rate $O(1/N^2)$ up to a logarithmic factor.
This is the best previously known rate,
and this paper provides a better rate.

In short, none of the existing first-order methods
achieve the optimal rate for the gradient inaccuracy even up to a constant,
and thus this paper focuses on optimizing the gradient efficiency
of first-order methods for smooth convex minimization
with %conditions
\eqref{eq:init,func}
and~\eqref{eq:init,func,N}.
Table~\ref{tab:rate} summarizes
the efficiency of first-order methods
and illustrates that the proposed \OGMG~attains 
the optimal worst-case gradient rate $O(1/N^2)$
%under the conditions~
with~\eqref{eq:init,func}
and~\eqref{eq:init,func,N}.

\begin{remark}
\label{rk:ogmg}
After the initial version of this paper was posted online~\cite{kim:18:otc},
the paper~\cite[Remark 2.1]{nesterov:20:pda}
constructed a simple method using \OGMG~that achieves $O(1/N^4)$
under the initial distance condition~\eqref{eq:init,dist}.
The method runs an accelerated method such as Nesterov's FGM and OGM
for the first half of the iterations
and then runs \OGMG~for the rest.
That approach (built upon the proposed \OGMG)
further closes the open problem of developing an optimal method
for decreasing the gradient,
under the initial distance condition~\eqref{eq:init,dist}. %and~\eqref{eq:init,func}.
%The best known gradient rate~\eqref{eq:epsgrad} under
%the initial distance condition~\eqref{eq:init,dist}
%(for both with or without knowing $\overline{R}$) 
%in Table~\ref{tab:rate} %and~\ref{tab:rateR}
%{\cblue can be now replaced by $O(1/N^4)$.} %and $O(L^2\overline{R}^2/N^4)$, respectively.
\end{remark}

\begin{table}[!h]
\centering
\caption{
Summary of the efficiency of first-order methods
discussed in Sect.~\ref{sec:eff}
\cite{nesterov:04,nemirovsky:92:ibc,nesterov:12:htm,carmon:19:lbf,taylor:18:ewc};
The rates of the proposed \OGMG~and a method in~\cite[Remark 2.1]{nesterov:20:pda}
using \OGMG~(see Remark~\ref{rk:ogmg}) are also presented.
$\widetilde{O}(\cdot)$ is a big-$O$ notation that ignores a logarithmic factor.
}
\label{tab:rate}
\begin{tabular}{|c|c|c|c|c|c|c|c|}
\hline
\multirow{2}{*}{Efficiency} %Inacurracy} 
	& Initial
	& GM 
	& \multicolumn{2}{c|}{Best known rate} 
	& OGM-G
	& \cite{nesterov:20:pda}
	& Optimal 
	\\
	\cline{4-5}

	& cond.
	& rate
	& w/o $R$ or $\overline{R}$
	& w/ $R$ or $\overline{R}$
	& rate
	& rate
	& rate
	\\
\hline
Cost func.~\eqref{eq:epsfunc}
	& \eqref{eq:init,dist} 
	& $O(1/N)$
	& \multicolumn{2}{c|}{$O(1/N^2)$}
	& $\cdot$
	& $\cdot$
	& $O(1/N^2)$
	\\
\hline
\multirow{3}{*}{Gradient~\eqref{eq:epsgrad}} 
	& \eqref{eq:init,dist}
	& $O(1/N^2)$
        & $O(1/N^3)$
	& $\widetilde{O}(1/N^4)$
	& $\cdot$
	& $O(1/N^4)$
        & $O(1/N^4)$
        \\
\cline{2-8}
	& \eqref{eq:init,func}
	& $O(1/N)$
	& $O(1/N)$
        & $\widetilde{O}(1/N^2)$
        & $O(1/N^2)$
	& $\cdot$
	& $O(1/N^2)$
        \\
\cline{2-8}
        & \eqref{eq:init,func,N}
        & $O(1/N)$
        & $\cdot$ %$O(1/N)$
        & $\cdot$ %$\widetilde{O}(1/N^2)$
        & $O(1/N^2)$
        & $\cdot$
        & $\cdot$ %{\cgreen $O(1/N^2)$}
	\\
\hline
\end{tabular}
\end{table}

As Table~\ref{tab:rate} demonstrates,
worst-case rates of any given method
and optimal worst-case rates 
depend dramatically on the initial condition.
In particular,
the worst-case gradient rates for~\eqref{eq:init,func}
tend to be slower than those for~\eqref{eq:init,dist}.
%
%In other words, the worst-case rates 
%measured on a given set of problem instances $(f,\x_0)$ 
%depends on the initial condition as well as the function class.
%
At first glance,
this situation might hinder one's interest 
on the initial function condition~\eqref{eq:init,func}
studied in this paper. 
However, one should also consider the constants $R$ and $\overline{R}$
%{\cblue as illustrated in} Table~\ref{tab:rateR}
for a fair comparison of the worst-case rates.
In particular, consider a problem instance $(f,\x_0)$ where $f\in\cF$
and $\Xs \neq \emptyset$.
%with $\x_*\in\Xs$.
Then, choose $R$ and $\overline{R}$ such that
\begin{align}
f(\x_0) - f_* = \frac{1}{2}LR^2
\quad\text{and}\quad 
||\x_0 - \x_*|| = \overline{R}
\end{align}
for some $\x_*\in\Xs$.
Using the inequality $f(\x_0) - f_* \le \frac{L}{2}||\x_0 - \x_*||^2$
due to the smoothness of $f$ in~\eqref{eq:Lsmooth},
we have the relationship: 
\begin{align}
R \le \overline{R}
.\end{align}
For any optimization method,
including GM %~\cite{nemirovsky:92:ibc,taylor:18:ewc} 
and OGM-G,
the ratio $\Frac{\overline{R}}{R}$ can be 
in the order of $N$
or beyond,
for a given $N$,
and should not be neglected.
Section~\ref{sec:worst,ogmg,dist} gives one such example.

\section{Performance Estimation Problem (PEP) 
for the Worst-Case Gradient Decrease}
\label{sec:pep}

This section studies PEP~\cite{drori:14:pof}
and its relaxations for the worst-case gradient analysis
under the condition~\eqref{eq:init,func}.

\subsection{\bf Exact PEP}

The papers~\cite{drori:14:pof,taylor:17:ssc} suggest that
for any given step coefficients $\h := \{h_{i,k}\}$ of a FSFOM,
total number of iterations $N$,
problem dimension $d$,
and constants $L$, $R$,
the exact worst-case gradient bound under~\eqref{eq:init,func} is given by
%\footnote{
%The condition $\Xs \neq \emptyset$ is included here 
%because the associated PEP proof assumes such condition. 
%Corollaries~\ref{cor:gm} and~\ref{cor:ogmg} derive worst-case bounds that do not require this condition
%for the alternative initial condition~\eqref{eq:init,func,N}.
%}
\begin{align}
\Bd_{\mathrm{P}}(\h,N,d,L,R) :=
&\max_{f\in\cF}
\max_{\x_0,\ldots,\x_N\in\reals^d} 
	\frac{1}{L^2R^2}||\nabla f(\x_N)||^2 
	\tag{P} 
	\label{eq:P} \\
&\st\; \begin{cases}
	\x_{i+1} = \x_i - \frac{1}{L}\sum_{k=0}^i h_{i+1,k}\nabla f(\x_k),
        \quad i=0,\ldots,N-1, & \\
	f(\x_0) - f_* \le \frac{1}{2}LR^2,
	\quad %f_* = \inf_{\x\in\reals^d} f(\x), \quad 
	\x_* \in \Xs \neq \emptyset, &
	\end{cases}
	\nonumber
\end{align}
where $||\nabla f(\x_N)||^2$ is multiplied by $\frac{1}{L^2R^2}$ for convenience in later analysis.
However, as noted in~\cite{drori:14:pof}, 
it is intractable to solve~\eqref{eq:P}
due to its infinite-dimensional function constraint.
Thus the next section employs relaxations introduced in~\cite{drori:14:pof}.

\subsection{\bf Relaxing PEP}

As suggested by~\cite{drori:14:pof,taylor:17:ssc},
to convert~\eqref{eq:P} into an equivalent finite-dimensional problem,
we replace the infinite-dimensional constraint $f\in\cF$ by
a finite set of inequalities satisfied by $f\in\cF$
%such $f$
\cite[Theorem~2.1.5]{nesterov:04}:
\begin{align}
\frac{1}{2L}||\nabla f(\x_i) - \nabla f(\x_j)||^2
\le f(\x_i) - f(\x_j) - \inprod{\nabla f(\x_j)}{\x_i - \x_j}
\label{eq:smooth}
\end{align}
on each pair $(i,j)$ for $i,j=0,\ldots,N,*$.
For simplicity in the proofs, we further narrow down the set\footnote{
We found that the set of constraints in~\eqref{eq:P1}
is sufficient for the exact worst-case gradient analysis
of GM and \OGMG~for~\eqref{eq:init,func}, %(and $d\ge N+1$),
as illustrated in later sections.
In other words, 
the resulting worst-case rates of GM and \OGMG~in this paper 
are tight with our specific choice of the set of inequalities. 
Note that
this relaxation choice in~\eqref{eq:P1}
differs from the choice in~\cite[Problem (G$'$)]{drori:14:pof}.} 
of inequalities~\eqref{eq:smooth},
specifically the pairs 
$\{(i-1,i)\;:\;i=1,\ldots,N\}$, 
$\{(N,i)\;:\;i=0,\ldots,N-1\}$
and $\{(N,*)\}$.\footnote{
\label{ft:nonempty}
The inequality~\eqref{eq:smooth}
for the pair $\{(N,*)\}$ simplifies to $\frac{1}{2L}||\nabla f(\x_N)||^2 \le f(\x_N) - f_*$
under the condition $\Xs \neq \emptyset$.
Such inequality is not used under the assumption~\eqref{eq:init,func,N} 
in Corollaries~\ref{cor:gm} and~\ref{cor:ogmg}.
}
This relaxation leads to
\begin{align}
\Bd_{\mathrm{P1}}(\h,N,d) :=
&\max_{\substack{\G\in\reals^{(N+1)\times d}, \\ \del\in\reals^{N+1}}}\;
        \Tr{\G^\top\u_N\u_N^\top\G}
	\tag{P1} 
	\label{eq:P1} \\
&\st\; \begin{cases}
	\Tr{\G^\top\A_{i-1,i}(\h)\G} \le \delta_{i-1} - \delta_i, \quad i=1,\ldots,N, & \\
	\Tr{\G^\top\B_{N,i}(\h)\G} \le \delta_N - \delta_i, \quad i=0,\ldots,N-1, & \\
	\Tr{\G^\top\C_N\G} \le \delta_N, \quad \delta_0 \le \frac{1}{2}, &
	\end{cases}
	\nonumber
\end{align}
where we define
\begin{align}
\begin{cases}
\g_i := \frac{1}{LR}\nabla f(\x_i), \quad i=0,\ldots,N, 
	\quad \G := [\g_0,\ldots,\g_N]^\top, & \\
\delta_i := \frac{1}{LR^2}(f(\x_i) - f_*), \quad i=0,\ldots,N, \quad
	\del := [\delta_0,\ldots,\delta_N]^\top, & \\
\u_i := [0,\ldots,0,\underbrace{1}_{\text{$(i+1)$th entry}},0,\ldots,0]^\top \in \reals^{N+1},
        \quad i=0,\ldots,N, &
\end{cases}
\label{eq:def1}
\end{align} 
and
\begin{align}
\begin{cases}
\A_{i-1,i}(\h) := \frac{1}{2}(\u_{i-1} - \u_i)(\u_{i-1} - \u_i)^\top
                + \frac{1}{2}\sum_{k=0}^{i-1}
		h_{i,k}(\u_i\u_k^\top + \u_k\u_i^\top), \quad i=1,\ldots,N, & \\
\B_{N,i}(\h) := \frac{1}{2}(\u_N - \u_i)(\u_N - \u_i)^\top
                - \frac{1}{2}\sum_{l=i+1}^N\sum_{k=0}^{l-1}
		h_{l,k}(\u_i\u_k^\top + \u_k\u_i^\top), & \\
\hspace{250pt} \quad i=0,\ldots,N-1, & \\
\C_N := \frac{1}{2}\u_N\u_N^\top.
\end{cases}
\label{eq:def2}
\end{align}

As in~\cite{taylor:17:ssc},
we further relax~\eqref{eq:P1} by introducing the Gram matrix $\Z := \G\G^\top$ as
\begin{align}
\Bd_{\mathrm{P2}}(\h,N,d) :=
&\max_{\substack{\Z\in\symm_+^{N+1}, \\ \del\in\reals^{N+1}}}\;
        \Tr{\u_N\u_N^\top\Z}
        \tag{P2} 
        \label{eq:P2} \\
&\st\; \begin{cases}
        \Tr{\A_{i-1,i}(\h)\Z} \le \delta_{i-1} - \delta_i, \quad i=1,\ldots,N, & \\
        \Tr{\B_{N,i}(\h)\Z} \le \delta_N - \delta_i, \quad i=0,\ldots,N-1, & \\
        \Tr{\C_N\Z} \le \delta_N, \quad
        \delta_0 \le \frac{1}{2}. &
        \end{cases}
        \nonumber
\end{align}
This problem has the following Lagrangian dual:
\begin{align}
\Bd_{\mathrm{D}}(\h,N) :=
&\min_{(\aa,\bb,c,e)\in\reals_+^{2N+2}}
        \frac{1}{2}e
	\tag{D}
	\label{eq:D} \\
&\st\; \begin{cases}
	\S(\h,\aa,\bb,c) \succeq\zero, \quad
	-a_1 + b_0 + e = 0, \quad a_N - \sum_{i=0}^{N-1}b_i - c = 0, & \\
        a_i - a_{i+1} + b_i = 0, \quad i=1,\ldots,N-1. & 
	\end{cases}
	\nonumber
\end{align}
where
\begin{align}
&\S(\h,\aa,\bb,c) 
        :=\; \sum_{i=1}^Na_i\A_{i-1,i}(\h) + \sum_{i=0}^{N-1}b_i\B_{N,i}(\h)
        + c\C_N(\h)
        - \u_N\u_N^\top 
        \label{eq:SS} \\
        =\;& \frac{1}{2}\sum_{i=1}^N a_i(\u_{i-1} - \u_i)(\u_{i-1} - \u_i)^\top
        + \frac{1}{2}\sum_{i=0}^{N-1} b_i(\u_N - \u_i)(\u_N - \u_i)^\top
        + \frac{1}{2}(c-2)\u_N\u_N^\top 
        \nonumber \\
        &+ \frac{1}{2}\sum_{i=1}^N\sum_{k=0}^{i-1} 
                a_i h_{i,k}
                (\u_i\u_k^\top + \u_k\u_i^\top)
        - \frac{1}{2}\sum_{i=0}^{N-1}\sum_{k=0}^{N-1}
                \paren{b_i\sum_{l=\max\left\{\substack{i+1,\\k+1}\right\}}^N h_{l,k}}
                (\u_i\u_k^\top + \u_k\u_i^\top)
        \nonumber
.\end{align}
For given $\h$ and $N$, a semidefinite programming (SDP) problem~\eqref{eq:D} 
can be solved numerically using an SDP solver
(\eg, CVX~\cite{cvxi,gb08}).
The next two sections analytically specify 
feasible points of~\eqref{eq:D} for GM and \OGMG,
which were numerically first identified to be solutions of~\eqref{eq:D}
for each method by the authors.
These feasible points provide 
the exact worst-case analytical gradient bounds for GM and \OGMG.

\section{Applying the Relaxed PEP to GM}
\label{sec:gm,ic2}

Inspired by the numerical solutions of~\eqref{eq:D} for GM
using CVX~\cite{cvxi,gb08},
we next specify a feasible point of~\eqref{eq:D} for GM.

\begin{lemma}
\label{lem:gm_feas_ic2}
For GM, \ie the FSFOM with $h_{i+1,k}$ having $1$ for $k=i$ and $0$ otherwise,
the following set of dual variables:
\begin{align}
\begin{cases}
a_i = \frac{2(N+i)}{(N-i+1)(2N+1)} = \frac{N+i}{N-i+1}e, \quad i = 1,\ldots,N, & \\
b_i = \begin{cases}
		\frac{2}{N(2N+1)} = \frac{1}{N}e, & i = 0, \\
		\frac{2}{(N-i)(N-i+1)}, & i = 1,\ldots,N-1,
	\end{cases} & \\
c = e = \frac{2}{2N+1}, &
\end{cases}
\label{eq:gm,feas}
\end{align}
is a feasible point of~\eqref{eq:D}.
\begin{proof}
Obviously,~\eqref{eq:gm,feas} %is in $\Lambda$~\eqref{eq:Lam},
satisfies the equality conditions of~\eqref{eq:D},
and the rest of proof shows the positive semidefinite condition of~\eqref{eq:D}.

For any $\h$ and $(\aa,\bb,c,e)\in\Lambda$, 
the $(i,j)$th entry of the symmetric matrix~\eqref{eq:SS} 
can be rewritten as
\begin{align}
&[2\S(\h,\aa,\bb,c)]_{ij} \label{eq:S} \\
=\;& \begin{cases}
        a_1 + b_0\paren{1 - 2\sum_{l=1}^N h_{l,0}}, 
                & i=0,\; j=i, \\
        a_i + a_{i+1} + b_i\paren{1 - 2\sum_{l=i+1}^N h_{l,i}}, 
                & i=1,\ldots,N-1,\; j=i, \\
        a_N + \sum_{l=0}^{N-1}b_l + c - 2 = 2(a_N - 1), 
                & i=N,\; j=i, \\
        a_i(h_{i,i-1} - 1) - b_i\sum_{l=i+1}^N h_{l,i-1} - b_{i-1}\sum_{l=i+1}^N h_{l,i}, 
                & i=1,\ldots,N-1,\; j=i-1, \\
        a_N(h_{N,N-1} - 1) - b_{N-1}, 
                & i=N,\; j=i-1, \\
        a_ih_{i,j} - b_i\sum_{l=i+1}^N h_{l,j} - b_j\sum_{l=i+1}^N h_{l,i}, 
		& i=2,\ldots,N-1, \\
		& j=0,\ldots,i-2, \\  
        a_Nh_{N,j} - b_j, 
                & i=N,\;j=0,\ldots,i-2.
\end{cases}
\nonumber
\end{align}
Substituting the step coefficients $\h$ for GM 
and the dual variables~\eqref{eq:gm,feas}
in~\eqref{eq:S} yields
\begin{align}
[2\S(\h,\aa,\bb,c)]_{ij}
&= \begin{cases}
        a_1 - b_0 = e, 
		& i=0,\; j=i, \\
        a_i + a_{i+1} - b_i = 2a_i, 
		& i=1,\ldots,N-1,\; j=i, \\
        2(a_N - 1), 
		& i=N,\; j=i, \\
        - b_j, 
		& i=1,\ldots,N,\;j=0,\ldots,i-1,
\end{cases}
\label{eq:S,gm}
\end{align}
The matrix~\eqref{eq:S,gm}
has nonnegative diagonal entries,
and thus showing the diagonal dominance of the matrix~\eqref{eq:S,gm}
implies its positive semidefiniteness.
%which then proves the positive semidefinite condition of~\eqref{eq:D} 
%with $\gamma=e$ in~\eqref{eq:gm,feas}.

A sum of absolute values of nondiagonal elements for each row is
\begingroup
\allowdisplaybreaks
\begin{align}
&\sum_{\substack{j=0 \\ j\neq i}}^N\abs{[2\S(\h,\aa,\bb,c)]_{ij}}
= \begin{cases}
Nb_0, & i=0, \\
b_0 + (N-1)b_1 
	%= \frac{2}{N(2N+1)} + \frac{2}{N} = \frac{4(N+1)}{N(2N+1)}, 
	& i=1, \\
\sum_{j=0}^{i-1}b_l + (N-i)b_i 
	%= \frac{2}{N(2N+1)} + \frac{2}{N-i+1} - \frac{2}{N} + \frac{2}{N-i+1}
	%= \frac{4(N+i)}{(N-i+1)(2N+1)}, 
	& i=2,\ldots,N-1, \\ 
\sum_{j=0}^{N-1}b_j
	%= \frac{2}{N(2N+1)} + 2 - \frac{2}{N}
	%= \frac{2(2N-1)}{2N+1}, 
	& i=N,
\end{cases} \\
=\; &\begin{cases}
\frac{2}{2N+1}, & i=0, \\
\frac{2}{N(2N+1)} + \frac{2}{N} = \frac{4(N+1)}{N(2N+1)}, & i=1, \\
\frac{2}{N(2N+1)} + \frac{2}{N-i+1} - \frac{2}{N} + \frac{2}{N-i+1}
        = \frac{4(N+i)}{(N-i+1)(2N+1)}, & i=2,\ldots,N-1, \\ 
\frac{2}{N(2N+1)} + 2 - \frac{2}{N}
        = \frac{2(2N-1)}{2N+1}, & i=N,
\end{cases}
\nonumber \\
=\; &\begin{cases}
e, & i=0, \\
\frac{2(N+i)}{(N-i+1)}e, & i=1,\ldots,N-1, \\
2(2Ne-1), & i=N,
\end{cases}
\nonumber
\end{align}
\endgroup
and this satisfies
$[2\S(\h,\aa,\bb,c)]_i
= \sum_{\substack{j=0 \\ j\neq i}}^N\abs{[2\S(\h,\aa,\bb,c)]_{ij}}$ for all $i$,
\ie, the matrix~\eqref{eq:S,gm} is diagonally dominant,
and this concludes the proof.
\qed
\end{proof}
\end{lemma}

The next theorem provides the worst-case convergence gradient bound of GM.
%which was first numerically identified in~\cite{taylor:18:ewc}.

\begin{theorem}
\label{thm:gm_bound_ic2}
Assume that $f\in\cF$, $\Xs \neq \emptyset$, 
and $f(\x_0) - f_* \le \frac{1}{2}LR^2$~\eqref{eq:init,func}.
Let $\x_0,\ldots,\x_N\in\reals^d$ be
generated by GM,
\ie, the FSFOM with $h_{i+1,k}$ having $1$ for $k=i$ and $0$ otherwise.
Then, for any $N\ge1$,
\begin{align}
||\nabla f(\x_N)||^2 \le \frac{L^2R^2}{2N+1}
\label{eq:gm_bound_ic2}
.\end{align}
\begin{proof}
Using Lemma~\ref{lem:gm_feas_ic2} for the step coefficients $\h$ of GM, we have
\begin{align*} 
||\nabla f(\x_N)||^2 
	\le L^2R^2\Bd_{\mathrm{D}}(\h,N) 
	\le L^2R^2\frac{1}{2N+1}
.\end{align*}
%Now apply~\eqref{eq:init,func}, \ie,
%let $f(\x_0) - f(\x_*) = \frac{1}{2}LR^2$.
\qed
\end{proof}
\end{theorem}

The PEP proof of Theorem~\ref{thm:gm_bound_ic2},
using Lemma~\ref{lem:gm_feas_ic2},
can be used to construct a conventional proof
that derives inequality~\eqref{eq:gm_bound_ic2}
by a weighted sum of the inequalities~\eqref{eq:smooth}.
Specifically, one can use a weighted sum of inequalities
using the dual variables $(\aa,\bb,c,e)$ in~\eqref{eq:gm,feas} as weights:
\begin{align}
\frac{1}{2L}||\nabla f(\x_{i-1}) - \nabla f(\x_i)||^2 
	\le f(\x_{i-1}) - f(\x_i) - \inprod{\nabla f(\x_i)}{\x_{i-1} - \x_i} \quad&:\quad a_i 
	\label{eq:ineqs} \\
\frac{1}{2L}||\nabla f(\x_N) - \nabla f(\x_i)||^2 
        \le f(\x_N) - f(\x_i) - \inprod{\nabla f(\x_i)}{\x_N - \x_i} \quad&:\quad b_i \nonumber \\
\frac{1}{2L}||\nabla f(\x_N)||^2 
        \le f(\x_N) - f_* \quad&:\quad c \nonumber \\
f(\x_0) - f_* \le \frac{1}{2}LR^2 \quad&:\quad e, \nonumber 
\end{align}
which simplifies to
\begin{align}
\frac{1}{L}||\nabla f(\x_N)||^2 
+ \sum_{i=1}^N\sum_{j=0}^{i-1}\frac{b_j}{2L}\left|\left|\nabla f(\x_i) - \nabla f(\x_j)\right|\right|^2 
\le \frac{LR^2}{2N+1}
\label{eq:gm,weighted}
,\end{align}
and this yields~\eqref{eq:gm_bound_ic2}.

We next show that the bound~\eqref{eq:gm_bound_ic2} is exact
by specifying a certain worst-case function.
This implies that the feasible point in~\eqref{eq:gm,feas}
is an optimal point of~\eqref{eq:D} for GM.

\begin{lemma}
\label{lem:worst,gm}
For the following Huber function in $\cF$ for all $d\ge1$:
\begin{align}
\phi(\x) = \begin{cases}
                \frac{LR}{\sqrt{2N+1}}||\x|| - \frac{LR^2}{2(2N+1)},
                        & ||\x|| \ge \frac{R}{\sqrt{2N+1}}, \\
                \frac{L}{2}||\x||^2, & \text{otherwise},
        \end{cases}
\end{align}
GM exactly achieves the bound~\eqref{eq:gm_bound_ic2}
with %an initial point 
$\x_0$ satisfying
$\phi(\x_0) - \phi_* = \frac{1}{2}LR^2$.
\begin{proof}
Starting from $\x_0 = \frac{N+1}{\sqrt{2N+1}}R\nnu$
that satisfies $\phi(\x_0) - \phi_* = \frac{1}{2}LR^2$~\eqref{eq:init,func}
for any unit-norm vector $\nnu$,
the iterates of GM are as follows:
\begin{align*}
\x_i = \x_0 - \frac{1}{L}\sum_{k=0}^{i-1} \nabla \phi(\x_k)
	= \paren{\frac{N+1}{\sqrt{2N+1}} - \frac{i}{\sqrt{2N+1}}}R\nnu,
	\quad i=0,\ldots,N,
\end{align*}
where all the iterates
stay in the affine region of the function $\phi(\x)$
with the same gradient $\nabla \phi(\x_i) = \frac{LR}{\sqrt{2N+1}}\nnu,\; i=0,\ldots,N$.
Therefore, after $N$ iterations of GM, we have
$ 
||\nabla \phi(\x_N)||^2 = \frac{L^2R^2}{2N+1}
,$
which concludes the proof.
\qed
\end{proof}
\end{lemma}

\begin{remark}
\label{rk:gm}
For $f\in\cF$,
and for some $\x_*\in\Xs$
and $||\x_0 - \x_*|| \le \overline{R}$~\eqref{eq:init,dist},
the $N$th iterate $\x_N$ of GM %, \ie,
%the FSFOM with $h_{i+1,k}$ having $1$ for $k=i$ and $0$ otherwise,
has the following exact
worst-case cost function bound~\cite[Theorems~1 and 2]{drori:14:pof}:
\begin{align}
f(\x_N) - f(\x_*) \le \frac{L\overline{R}^2}{2(2N+1)}
\label{eq:gm_cost}
,\end{align}
where this exact upper bound is
equivalent to the exact worst-case gradient bound~\eqref{eq:gm_bound_ic2} of GM
up to a constant $\frac{\overline{R}^2}{2LR^2}$.
A similar relationship appears in
\cite[Table 3]{taylor:18:ewc}
for nonsmooth composite convex minimization.
\end{remark}

The preceding results in this section assume 
that there is a finite minimizer. %$\x_*\in\Xs$.
There are applications that do not have a finite minimizer $\x_*\in\Xs$,
\eg, an unregularized logistic regression of an overparameterized model for separable datasets
\cite{nacson:19:cog,soudry:18:tib}.
The following corollary extends the analysis to such cases.

\begin{corollary}
\label{cor:gm}
For $f\in\cF$,
let $\x_0,\ldots,\x_N\in\reals^d$ be
generated by GM.
Assume that $f(\x_0) - f(\x_N) \le \frac{1}{2}LR_N^2$~\eqref{eq:init,func,N}.
Then, for any $N\ge1$,
\begin{align}
||\nabla f(\x_N)||^2 \le \frac{L^2R_N^2}{2N}
\label{eq:gm,N}
.\end{align}
\begin{proof}
Equation~\eqref{eq:gm,weighted} 
consists of a weighted sum of the third and fourth inequalities of~\eqref{eq:ineqs},
scaled by $c=e=\frac{2}{2N+1}$ in~\eqref{eq:gm,feas}:
\begin{align*}
\frac{c}{2L}||\nabla f(\x_N)||^2 + c(f(\x_0) - f(\x_N)) \le \frac{c}{2}LR^2
.\end{align*}
The third inequality of~\eqref{eq:ineqs} 
assumes $\Xs \neq \emptyset$ (see footnote~\ref{ft:nonempty}),
so we derive a bound without the above inequality.
Replacing the above inequality in the weighted summation 
for deriving~\eqref{eq:gm,weighted} by
\eqref{eq:init,func,N} 
%(that does not involve $\x_*$)} 
scaled by $c$,
$
c(f(\x_0) - f(\x_N)) \le \frac{c}{2}LR_N^2
$,
yields
\begin{align*}
\frac{1}{L}\left(1 - \frac{c}{2}\right)||\nabla f(\x_N)||^2
+ \sum_{i=1}^N\sum_{j=0}^{i-1}\frac{b_j}{2L}\left|\left|\nabla f(\x_i) - \nabla f(\x_j)\right|\right|^2
\le \frac{LR_N^2}{2N+1}
,\end{align*}
which concludes the proof.
\qed
\end{proof}
\end{corollary}

\section{Optimizing FSFOM Using the Relaxed PEP}
\label{sec:optpep}

This section optimizes the step coefficients of FSFOM 
using the relaxed PEP~\eqref{eq:D}
to develop an efficient first-order method
for decreasing the gradient of smooth convex functions.

\subsection{\bf Numerically Optimizing FSFOM Using the Relaxed PEP}

To optimize the step coefficients of $\h$ of FSFOM for each given $N$, 
we are interested in solving
\begin{align}
\tilde{\h} := \argmin{\h} \Bd_{\mathrm{D}}(\h,N)
\tag{HD}
\label{eq:HD}
,\end{align}
which is nonconvex.
However, the problem~\eqref{eq:HD} is bi-convex
over $\h$ and $(\aa,\bb,c,e,\gamma)$,
so for each given $N$
we numerically solved~\eqref{eq:HD} by an alternating minimization approach
using CVX~\cite{cvxi,gb08}.
Inspired by those numerical solutions,
the next section specifies a feasible point of~\eqref{eq:HD}.

\subsection{\bf A Feasible Point of the Relaxed PEP}

The following lemma specifies a feasible point of~\eqref{eq:HD}.

\begin{lemma}
\label{lem:ogmg_feas_ic2}
The following step coefficients of FSFOM:
\begin{align}
\tth_{i+1,k} 
	&= \begin{cases}
        \frac{\ttheta_{k+1} - 1}{\ttheta_k}\tth_{i+1,k+1}, & k=0,\ldots,i-2, \\
        \frac{\ttheta_{k+1} - 1}{\ttheta_k}(\tth_{i+1,i} - 1), & k=i-1, \\
        1 + \frac{2\ttheta_{i+1} - 1}{\ttheta_i}, & k=i,
\end{cases}
\label{eq:ogmg,h}
\end{align}
and the following set of dual variables:
\begin{align}
%\begin{cases}
a_i = \frac{1}{\ttheta_i^2}, \quad i = 1,\ldots,N, \quad 
b_i = \frac{1}{\ttheta_i\ttheta_{i+1}^2}, \quad i = 0,\ldots,N-1, \quad
c = e = \frac{2}{\ttheta_0^2}, 
%\end{cases}
\label{eq:ogmg,feas}
\end{align}
constitute a feasible point of~\eqref{eq:HD}
for the parameters: 
\begin{align}
\ttheta_i = \begin{cases}
                \frac{1 + \sqrt{1 + 8\ttheta_{i+1}^2}}{2}, & i=0, \\
                \frac{1 + \sqrt{1 + 4\ttheta_{i+1}^2}}{2}, & i=1,\ldots,N-1, \\
                1, & i = N.
        \end{cases} 
\label{eq:ttheta}
\end{align}
\begin{proof}
The appendix %~\ref{appx:ogmg,h} 
first derives 
properties of the step coefficients $\tilde{\h} = \{\tth_{i,k}\}$~\eqref{eq:ogmg,h}
that are used in the proof:
\begin{align}
\tth_{i,j} &= \frac{\ttheta_i^2(2\ttheta_i-1)}{\ttheta_j\ttheta_{j+1}^2}, 
	\quad i=2,\ldots,N,\;j=0,\ldots,i-2, 
	\label{eq:ogmg,h,sum1} \\
\sum_{l=i+1}^N \tth_{l,j} 
        &= \begin{cases}
		\frac{1}{2}(\ttheta_0+1), 
			& i=0,\;j=i, \\
		\ttheta_i,
			& i=1,\ldots,N-1,\;j=i, \\
		\frac{\ttheta_{i+1}^4}{\ttheta_j\ttheta_{j+1}^2}, 
                	& i=1,\ldots,N-1,\;j=0,\ldots,i-1. \\
	\end{cases}
	\label{eq:ogmg,h,sum2}
\end{align}
By definition of $\ttheta_i$~\eqref{eq:ttheta},
we also have
\begin{align}
\ttheta_i^2 &= \begin{cases}
                \ttheta_i + 2\ttheta_{i+1}^2, & i=0, \\
                \ttheta_i + \ttheta_{i+1}^2, & i=1,\ldots,N-1.
                \end{cases}
        \label{eq:ttheta,rule}
\end{align}

Obviously,~\eqref{eq:ogmg,feas} %is in $\Lambda$~\eqref{eq:Lam},
satisfies the equality conditions of~\eqref{eq:D},
and the rest of proof shows the positive semidefinite condition of~\eqref{eq:D}.
Substituting the step coefficients $\tilde{\h}$~\eqref{eq:ogmg,h} 
and the dual variables~\eqref{eq:ogmg,feas}
with their properties~\eqref{eq:ogmg,h,sum1},~\eqref{eq:ogmg,h,sum2}
and~\eqref{eq:ttheta,rule}
in~\eqref{eq:S} yields
\begin{align*}
&\;[2\S(\h,\aa,\bb,c)]_{ij} \\
=& \begin{cases}
        \frac{1}{\ttheta_1^2} + \frac{1}{\ttheta_0\ttheta_1^2}(1-(\ttheta_0+1)),
		\quad i=0,\; j=i, & \\
        \frac{1}{\ttheta_i^2} + \frac{1}{\ttheta_{i+1}^2} 
		+ \frac{1}{\ttheta_i\ttheta_{i+1}^2}
		\paren{1 - 2\ttheta_i} 
		= \frac{\ttheta_{i+1}^2 + \ttheta_i - \ttheta_i^2}{\ttheta_i^2\ttheta_{i+1}^2}, 
		\quad i=1,\ldots,N-1,\; j=i, & \\
        2\paren{\frac{1}{\ttheta_N^2} - 1}, 
		\quad i=N,\; j=i, & \\
        \frac{1}{\ttheta_i^2}\frac{2\ttheta_i-1}{\ttheta_{i-1}} 
		- \frac{1}{\ttheta_i\ttheta_{i+1}^2}\frac{\ttheta_{i+1}^4}{\ttheta_{i-1}\ttheta_i^2}
		- \frac{1}{\ttheta_{i-1}\ttheta_i^2}\ttheta_i
		= \frac{(2\ttheta_i-1)\ttheta_i-\ttheta_{i+1}^2 - \ttheta_i^2}{\ttheta_{i-1}\ttheta_i^3},
		\; i=1,\ldots,N-1,\; j=i-1, & \\
	\frac{1}{\ttheta_N^2}\frac{2\ttheta_N-1}{\ttheta_{N-1}}	
		- \frac{1}{\ttheta_{N-1}\ttheta_N^2},
		\quad i=N,\; j=i-1, &\\
        \frac{1}{\ttheta_i^2}\frac{\ttheta_i^2(2\ttheta_i-1)}{\ttheta_j\ttheta_{j+1}^2} 
		- \frac{1}{\ttheta_i\ttheta_{i+1}^2}\frac{\ttheta_{i+1}^4}{\ttheta_j\ttheta_{j+1}^2} 
		- \frac{1}{\ttheta_j\ttheta_{j+1}^2}\ttheta_i
		= \frac{(2\ttheta_i-1)\ttheta_i - (\ttheta_i-1)^2\ttheta_i^2 - \ttheta_i^2}{\ttheta_j\ttheta_{j+1}^2\ttheta_i},
		\;\;\; i=2,\ldots,N-1, \\
		\hspace{253pt} j=0,\ldots,i-2, & \\
        \frac{1}{\ttheta_N^2}\frac{1}{\ttheta_j\ttheta_{j+1}^2}
		- \frac{1}{\ttheta_j\ttheta_{j+1}^2}, 
		\quad i=N,\; j=0,\ldots,i-2, &
\end{cases} \\
=&\; \zero,
%=& \begin{cases}
%	\frac{(2\ttheta_i-1)\ttheta_i-\ttheta_{i+1}^2 - \ttheta_i^2}{\ttheta_{i-1}\ttheta_i^3} 
%	= 0, & i=1,\ldots,N-1,\; j=i-1, \\
%	\frac{(2\ttheta_i-1)\ttheta_i - (\ttheta_i-1)^2\ttheta_i^2 - \ttheta_i^2}{\ttheta_j\ttheta_{j+1}^2\ttheta_i}
%	= 0, & i=2,\ldots,N-1,\; j=0,\ldots,i-2, \\
%	0, & i=0,\ldots,N,\;j=i, \\
%	0, & i=N,\;j=0,\ldots,i-1,
%\end{cases}	
\end{align*}
%This equality shows 
%the positive semidefinite condition of~\eqref{eq:D} (and~\eqref{eq:HD}),
%%%%%%using $\gamma = e$,
which concludes the proof.
\qed
\end{proof}
\end{lemma}

The next theorem provides the worst-case convergence gradient bound of 
FSFOM with step coefficients~\eqref{eq:ogmg,h}.

\begin{theorem}
\label{thm:ogmg,rate}
Assume that $f\in\cF$, $\Xs \neq \emptyset$, 
and $f(\x_0) - f_* \le \frac{1}{2}LR^2$~\eqref{eq:init,func}.
Let $\x_0,\ldots,\x_N\in\reals^d$ be
generated by FSFOM with step coefficients~\eqref{eq:ogmg,h}.
Then, for any $N\ge1$,
\begin{align}
||\nabla f(\x_N)||^2 \le \frac{L^2R^2}{\ttheta_0^2}
\le \frac{2L^2R^2}{(N+1)^2}
\label{eq:ogmg,rate}
.\end{align}
\begin{proof}
Using Lemma~\ref{lem:ogmg_feas_ic2}, we have
$ 
||\nabla f(\x_N)||^2 
        \le L^2R^2\Bd_{\mathrm{D}}(\h,N) 
        \le L^2R^2\frac{1}{\ttheta_0^2}
.$ 
We can easily show that $\ttheta_i$~\eqref{eq:ttheta}
satisfies $\ttheta_i \ge \frac{N-i+2}{2}$
for $i=1,\ldots,N$ by induction,
and this then yields
$\ttheta_0 \ge \frac{N+1}{\sqrt{2}}$,
%Now apply~\eqref{eq:init,func}, \ie,
%let $f(\x_0) - f(\x_*) = \frac{1}{2}LR^2$,
which concludes the proof.
\qed
\end{proof}
\end{theorem}

Similar to~\eqref{eq:gm,weighted},
the PEP proof of Theorem~\ref{thm:ogmg,rate},
using Lemma~\ref{lem:ogmg_feas_ic2},
can be used to construct a conventional proof
by a weighted sum of inequalities~\eqref{eq:ineqs}
using the dual variables $(\aa,\bb,c,e)$ in~\eqref{eq:ogmg,feas} 
as weights. This weighted sum leads to
\begin{align}
\frac{1}{L}||\nabla f(\x_N)||^2 \le \frac{LR^2}{\ttheta_0^2}
\label{eq:ogmg,weighted}
\end{align}
and yields~\eqref{eq:ogmg,rate}.

The bound~\eqref{eq:ogmg,rate} of FSFOM with~\eqref{eq:ogmg,h}
is optimal up to a constant
because Nemirovsky shows in~\cite{nemirovsky:92:ibc}
that the worst-case rate for the gradient decrease
of large-dimensional convex \emph{quadratic} function is $O(1/N^2)$
under~\eqref{eq:init,func}.
This result fills in Table~\ref{tab:rate}, %and~\ref{tab:rateR},
improving upon best known rates.

The following corollary examines the rate of FSFOM with~\eqref{eq:ogmg,h}
for cases where a finite minimizer might not exist.

\begin{corollary}
\label{cor:ogmg}
For $f\in\cF$,
let $\x_0,\ldots,\x_N\in\reals^d$ be
generated by FSFOM with step coefficients~\eqref{eq:ogmg,h}.
Assume that $f(\x_0) - f(\x_N) \le \frac{1}{2}LR^2$~\eqref{eq:init,func,N}.
Then, for any $N\ge1$,
\begin{align}
||\nabla f(\x_N)||^2 \le \frac{L^2R_N^2}{\ttheta_0^2 - 1}
\end{align}
\begin{proof}
Equation~\eqref{eq:ogmg,weighted} consists of
a weighted sum of the third and fourth inequalties of~\eqref{eq:ineqs},
scaled by $c=e=\frac{2}{\ttheta_0^2}$ in~\eqref{eq:ogmg,feas}:
\begin{align}
\frac{c}{2L}||\nabla f(\x_N)||^2 + c(f(\x_0) - f(\x_N)) \le \frac{c}{2}LR^2
.\end{align}
The third inequality of~\eqref{eq:ineqs} 
assumes $\Xs \neq \emptyset$ (see footnote~\ref{ft:nonempty}),
so we derive a bound without the above inequality.
Replacing the above inequality in the weighted summation 
for deriving~\eqref{eq:ogmg,weighted}
by~\eqref{eq:init,func,N} 
%(that does not involve $\x_*$)} 
scaled by $c$,
$c(f(\x_0) - f(\x_N)) \le \frac{c}{2}LR_N^2$, yields
\begin{align*}
\frac{1}{L}\left(1 - \frac{c}{2}\right)||\nabla f(\x_N)||^2
\le \frac{LR_N^2}{\ttheta_0^2}
,\end{align*}
which concludes the proof.
\qed
\end{proof}
\end{corollary}

The per-iteration computational complexity of the FSFOM
with~\eqref{eq:ogmg,h} 
would be expensive 
if implemented directly via~\eqref{eq:fsfom},
compared to GM, FGM and OGM,
so the next section provides an efficient form.

\subsection{\bf An Efficient Form of the Proposed Optimized Method: \OGMG}

This section develops an efficient form of FSFOM 
with the step coefficients~\eqref{eq:ogmg,h},
named \OGMG.
This form is similar to that of OGM~\cite{kim:16:ofo},
which is further studied in Sect.~\ref{sec:relatedOGM}.

\fbox{
\begin{minipage}[t]{0.85\linewidth}
\vspace{-10pt}
\begin{flalign*}
&\quad \text{\bf OGM-G} & \\
&\qquad \text{Input: } f\in\cF,\; \x_0=\y_0\in\reals^d,\; N\ge1. & \\
&\qquad \ttheta_i = \begin{cases}
		\frac{1 + \sqrt{1 + 8\ttheta_{i+1}^2}}{2}, & i=0, \\
		\frac{1 + \sqrt{1 + 4\ttheta_{i+1}^2}}{2}, & i=1,\ldots,N-1, \\
		1, & i = N,
	\end{cases} & \\
&\qquad \text{For } i = 0,\ldots,N-1, & \\
&\qquad \qquad \y_{i+1} = \x_i - \frac{1}{L}\nabla f(\x_i), & \\
&\qquad \qquad \x_{i+1} = \y_{i+1}
                + \frac{(\ttheta_i-1)(2\ttheta_{i+1} - 1)}
			{\ttheta_i(2\ttheta_i - 1)}(\y_{i+1} - \y_i)
                + \frac{2\ttheta_{i+1}-1}{2\ttheta_i-1}(\y_{i+1} - \x_i). &
%&\qquad \qquad \x_{i+1} = \y_{i+1}
%                + \frac{(2\ttheta_{i+1} - 1)\ttheta_{i-1}}
%                        {\ttheta_i(2\ttheta_i - 1)}
%		\left[\frac{\ttheta_i - 1}{\ttheta_{i-1}}(\y_{i+1} - \y_i)
%                + \frac{\ttheta_i}{\ttheta_{i-1}}(\y_{i+1} - \x_i)\right]. & %\\
%&\qquad \qquad \x_{i+1} = \y_{i+1}
%                + \frac{(2\ttheta_{i+1} - 1)\ttheta_{i+1}}
%                        {\ttheta_i(2\ttheta_i - 1)}
%                \left[\frac{\ttheta_i - 1}{\ttheta_{i+1}}(\y_{i+1} - \y_i)
%                + \frac{\ttheta_i}{\ttheta_{i+1}}(\y_{i+1} - \x_i)\right]. &
\end{flalign*}
\end{minipage}
} \vspace{5pt}

\begin{proposition}
\label{prop:ogmg}
The sequence $\{\x_0,\ldots,\x_N\}$ generated by FSFOM with~\eqref{eq:ogmg,h}
is identical to the corresponding sequence generated by \OGMG.
\begin{proof}
We first show that the step coefficients $\{\tth_{i+1,k}\}$~\eqref{eq:ogmg,h} 
are equivalent to
\begin{align}
\tth'_{i+1,k} 
= \begin{cases}
        \frac{(\ttheta_i - 1)(2\ttheta_{i+1} - 1)}{\ttheta_i(2\ttheta_i - 1)}\tth'_{i,k}, 
                & k=0,\ldots,i-2, \\
        \frac{(\ttheta_i - 1)(2\ttheta_{i+1} - 1)}{\ttheta_i(2\ttheta_i - 1)}(\tth'_{i,i-1} - 1), 
                & k=i-1, \\
        1 + \frac{2\ttheta_{i+1} - 1}{\ttheta_i}, & k=i.
\end{cases} 
\label{eq:ogmg,hh}
\end{align}
Obviously, $\tth_{i+1,i} = \tth'_{i+1,i},\;i=0,\ldots,N-1$, and we have
\begin{align*}
\tth_{i+1,i-1} &= \frac{\ttheta_i-1}{\ttheta_{i-1}}(\tth_{i+1,i} - 1)
        = \frac{(\ttheta_i-1)(2\ttheta_{i+1}-1)}{\ttheta_{i-1}\ttheta_i}
        = \frac{(\ttheta_i-1)(2\ttheta_{i+1}-1)}{\ttheta_i(2\ttheta_i-1)}
                \frac{2\ttheta_i-1}{\ttheta_{i-1}} \\
        &= \frac{(\ttheta_i-1)(2\ttheta_{i+1}-1)}{\ttheta_i(2\ttheta_i-1)}
		(\tth'_{i,i-1}-1)
	= \tth'_{i+1,i-1}
\end{align*}
for $i=1,\ldots,N-1$.

We next use induction by assuming
$\tth_{i+1,k}=\tth'_{i+1,k}$ for $i=0,\ldots,n-1,\;k=0,\ldots,i$.
We then have
\begin{align*}
\tth_{n+1,k} &= \frac{\ttheta_{k+1}-1}{\ttheta_k}\tth_{n+1,k+1}
        = \paren{\prod_{j=k}^{n-1}\frac{\ttheta_{l+1}-1}{\ttheta_l}}(\tth_{n+1,n}-1) \\
        &= \paren{\prod_{j=k}^{n-2}\frac{\ttheta_{l+1}-1}{\ttheta_l}}(\tth_{n,n-1}-1)
                \frac{\ttheta_n-1}{\ttheta_{n-1}}\frac{\tth_{n+1,n}-1}{\tth_{n,n-1}-1} \\
        &= \tth_{n,k}\frac{\ttheta_n-1}{\ttheta_{n-1}}
		\frac{(2\ttheta_{n+1} - 1)\ttheta_{n-1}}{\ttheta_n(2\ttheta_n-1)}
        = \frac{(\ttheta_n-1)(2\ttheta_{n+1}-1)}{\ttheta_n(2\ttheta_n-1)}\tth'_{n,k}
	= \tth'_{n+1,k}
\end{align*}
for $k=0,\ldots,n-2$,
where the fourth equality uses the definition of $\tth_{n,k}$.
This proves the first claim that the step coefficients 
$\{\tilde{h}_{i+1,k}\}$~\eqref{eq:ogmg,h}
and $\{\tilde{h}'_{i+1,k}\}$~\eqref{eq:ogmg,hh}
are equivalent.

We finally use induction to show the equivalence between
the generated sequences of
FSFOM with~\eqref{eq:ogmg,hh} and \OGMG.
For clarity, we use the notation $\x_0',\ldots,\x_N'$ and $\y_0',\ldots,\y_N'$
for \OGMG.
Obviously, $\x_0 = \x_0'$, and we have
\begingroup
\allowdisplaybreaks
\begin{align*}
\x_1 &= \x_0 - \frac{1}{L}\tth'_{1,0}\nabla f(\x_0)
	= \x_0 - \frac{1}{L}\paren{1 + \frac{2\ttheta_1-1}{\ttheta_0}}\nabla f(\x_0) \\ 
	&= \y_1' - \frac{1}{L}\frac{(2\ttheta_0-1)(2\ttheta_1-1)}
			{\ttheta_0(2\ttheta_0 - 1)}\nabla f(\x_0') \\
	&= \y_1' - \frac{1}{L}\paren{\frac{(\ttheta_0-1)(2\ttheta_1-1)}{\ttheta_0(2\ttheta_0-1)}
		+ \frac{2\ttheta_1-1}{2\ttheta_0-1}}\nabla f(\x_0')
	= \x_1'.
\end{align*}
\endgroup
Assuming $\x_i=\x_i'$ for $i=0,\ldots,n$, we have
\begingroup
\allowdisplaybreaks
\begin{align*}
&\x_{n+1} = \x_n - \frac{1}{L}\tth'_{n+1,n}\nabla f(\x_n)
		- \frac{1}{L}\tth'_{n+1,n-1}\nabla f(\x_{n-1})
		- \frac{1}{L}\sum_{k=0}^{n-2}\tth'_{n+1,k}\nabla f(\x_k) \\
	=\;& \x_n - \frac{1}{L}\paren{1 + \frac{2\ttheta_{n+1}-1}{\ttheta_n}}\nabla f(\x_n)
		- \frac{1}{L}\frac{(\ttheta_n-1)(2\ttheta_{n+1} - 1)}
			{\ttheta_n(2\ttheta_n - 1)}(\tth_{n,n-1} - 1)\nabla f(\x_{n-1}) \\
	&\quad
		- \frac{1}{L}\frac{(\ttheta_n-1)(2\ttheta_{n+1} - 1)}
                        {\ttheta_n(2\ttheta_n - 1)}\sum_{k=0}^{n-2}\tth_{n,k}\nabla f(\x_k) \\
	=\;& \x_n - \frac{1}{L}\paren{1 + \frac{2\ttheta_{n+1} - 1}
			{2\ttheta_n - 1}}\nabla f(\x_n) \\
	&\quad	+ \frac{(\ttheta_n-1)(2\ttheta_{n+1} - 1)}{\ttheta_n(2\ttheta_n - 1)}
			\paren{-\frac{1}{L}\nabla f(\x_n) + \frac{1}{L}\nabla f(\x_{n-1})
			- \frac{1}{L}\sum_{k=0}^{n-1}\tth_{n,k}\nabla f(\x_k)} \\
	=\;& \y_{n+1}' + \frac{(\ttheta_n-1)(2\ttheta_{n+1} - 1)}
			{\ttheta_n(2\ttheta_n - 1)}(\y_{n+1}' - \y_n')
		+ \frac{2\ttheta_{n+1} - 1}{2\ttheta_n - 1}(\y_{n+1}' - \x_n')
	= \x_{n+1}'.
\end{align*}
\endgroup
\qed
\end{proof}
\end{proposition}

\subsection{\bf Two Worst-Case Iterative Behaviors of \OGMG}
\label{sec:optpep,worst}

This section specifies two worst-case problem instances for \OGMG,
associated with Huber and quadratic functions respectively,
that make the bound~\eqref{eq:ogmg,rate} exact.
These examples imply that the feasible point in~\eqref{eq:ogmg,feas}
is an optimal point of~\eqref{eq:D} for \OGMG.

\begin{lemma}
\label{lem:worst}
For the following Huber and quadratic functions in $\cF$: %for all $d\ge1$:
\begin{align}
\phi_1(\x) = \begin{cases}
		\frac{LR}{\ttheta_0}||\x|| - \frac{LR^2}{2\ttheta_0^2}, 
			& ||\x|| \ge \frac{R}{\ttheta_0}, \\
		\frac{L}{2}||\x||^2, & \text{otherwise},
	\end{cases}
\quad
\text{and} 
\quad
\phi_2(\x) = \frac{L}{2}||\x||^2,
\label{eq:worst,two}
\end{align}
for all $d\ge1$,
\OGMG~exactly achieves the bound~\eqref{eq:ogmg,rate}
with an initial point 
$\x_0$ satisfying 
$\phi_1(\x_0) - \phi_{1,*} = \phi_2(\x_0) - \phi_{2,*} = \frac{1}{2}LR^2$.
\begin{proof}
We first consider $\phi_1(\x)$. 
Starting from an initial point $\x_0 = \frac{\ttheta_0^2+1}{2\ttheta_0}R\nnu$
that satisfies $\phi_1(\x_0) - \phi_{1,*} = \frac{1}{2}LR^2$~\eqref{eq:init,func}
for any unit-norm vector $\nnu$, 
we have
\begin{align*}
\x_N = \x_0 - \frac{1}{L}\sum_{j=1}^N\sum_{k=0}^{j-1}\tth_{j,k}\nabla f(\x_k)
	= \paren{\frac{\ttheta_0^2+1}{2\ttheta_0} 
		- \frac{\ttheta_0^2-1}{2\ttheta_0}}R\nnu
,\end{align*}
since
\begin{align*}
\sum_{j=1}^N\sum_{k=0}^{j-1}\tth_{j,k}
= \frac{1}{2}(\ttheta_0 + 1) + \sum_{j=1}^{N-1}\ttheta_j
= \frac{1}{2}(\ttheta_0 + 1 + 2\ttheta_1^2 - 2)
= \frac{1}{2}(\ttheta_0^2-1)
\end{align*}
that uses~\eqref{eq:ogmg,h,sum2} and~\eqref{eq:ttheta,rule}.
Here, all the iterates stay in the affine region of the function $\phi_1(\x)$
with the same gradient $\nabla\phi_1(\x) = \frac{LR}{\ttheta_0}\nnu,\; i=0,\ldots,N$.
Therefore, after $N$ iterations of OGM-G, we have
$ 
||\nabla\phi_1(\x_N)||^2 = \frac{L^2R^2}{\ttheta_0^2}
.$ 

We next consider $\phi_2(\x)$.
Starting from an initial point $\x_0 = R\nnu$ 
that satisfies $\phi_2(\x_0) - \phi_{2,*} = \frac{1}{2}LR^2$~\eqref{eq:init,func}
for any unit-norm vector $\nnu$, 
we have
\begin{align*}
\x_1 &= -\frac{1}{L}\frac{2\ttheta_1 - 1}{\ttheta_0}\nabla f(\x_0) 
	= -\frac{2\ttheta_1 - 1}{\ttheta_0}\x_0
,\end{align*}
and we have
\begin{align*}
\x_{i+1} &= - \frac{1}{L}\frac{2\ttheta_{i+1}-1}{2\ttheta_i-1}\nabla f(\x_i)
	= - \frac{2\ttheta_{i+1}-1}{2\ttheta_i-1}\x_i 
	= (-1)^i\frac{2\ttheta_{i+1}-1}{2\ttheta_1-1}\x_1,
	\;\; i=1,\ldots,N-1,
\end{align*}
using $\y_i = \zero,\;i=1,\ldots,N$.
Therefore, %after $N$ iterations of OGM-G, 
we have
$ 
||\nabla \phi_2(\x_N)||^2 = L^2||\x_N||^2 = \frac{L^2R^2}{\ttheta_0^2}
.$ 
\qed
\end{proof}
\end{lemma}

The iterates of OGM-G for the Huber worst-case function $\phi_1$
stay in one side of the affine region of the function,
while those for the quadratic worst-case function $\phi_2$
always overshoot the optimum.
These are extreme cases, and 
it is notable that some other first-order methods
also have two such worst-case iterative behaviors.
Specifically,
in~\cite{taylor:17:ssc,kim:17:otc}, first-order methods 
that have such two types of worst-case iterative behaviors
in Lemma~\ref{lem:worst},
associated with Huber and quadratic functions, respectively,
were found to have an optimal worst-case bound 
among a certain subset of first-order methods.
%a gradient method with the optimal constant step size and OGM 
%have such two types of worst-case behaviors,
%and have an optimal worst-case bound
%among fixed-step gradient methods and all first-order methods 
%(under a large-dimensional condition) respectively.}
This leads us to conjecture that the exact worst-case bound~\eqref{eq:ogmg,rate} 
of OGM-G may be optimal,
but proving it remains an open problem.

\subsection{\bf Worst-Case Rate Behaviors of OGM-G under Initial Distance Condition}
\label{sec:worst,ogmg,dist}

This section further studies the worst-case rate behaviors of OGM-G
under initial distance condition~\eqref{eq:init,dist}.
Table~\ref{tab:ratePESTO}
presents exact numerical worst-case rates of OGM-G
(under a large-dimensional condition),
using the performance estimation toolbox, named PESTO\footnote{
In PESTO toolbox~\cite{taylor:17:pet},
we used the SDP solver SeDuMi~\cite{sturm:99:us1} 
interfaced through Yalmip~\cite{lofberg:04:yat}.
The \OGMG~method is implemented in the PESTO toolbox.
}~\cite{taylor:17:pet},
based on PEP~\cite{drori:14:pof,taylor:17:ssc}.

\begin{table}[!h]
\centering
\caption{
Exact values of the reciprocals of
the worst-case cost function inaccuracy
$\paren{\frac{L\overline{R}^2}{f(\x_N)-f(\x_*)}}$
in~\eqref{eq:epsfunc}
and the worst-case gradient inaccuracy
$\paren{\frac{L^2\overline{R}^2}{||\nabla f(\x_N)||^2},
\frac{L^2R^2}{||\nabla f(\x_N)||^2},
\text{ or } \frac{L^2R_N^2}{||\nabla f(\x_N)||^2}}$
in~\eqref{eq:epsgrad}
of OGM-G
under one of the conditions~\eqref{eq:init,dist},~\eqref{eq:init,func}
or~\eqref{eq:init,func,N}.
}
\label{tab:ratePESTO}
\begin{tabular}{|c|c|c|c|c|c|c|c|c|c|} %c|}
\hline
\multirow{2}{*}{OGM-G Efficiency} %Inacurracy} 
        & Initial
        & \multicolumn{8}{c|}{Number of iterations}
        %& Empirical
        \\
        \cline{3-10}
        & cond.
        & 1
        & 2
        & 4
        & 10
        & 20
        & 30
        & 40
        & 50
        %& rate
        \\
\hline
Cost func.~\eqref{eq:epsfunc}
        & \eqref{eq:init,dist}
        & 8.0
        & 10.0
        & 9.7
        & 8.9
        & 8.5
        & 8.3
        & 8.3
        & 8.2
        %& $O\paren{\frac{L\overline{R}^2}{?}}$
        \\
\hline
\multirow{3}{*}{Gradient~\eqref{eq:epsgrad}}
        & \eqref{eq:init,dist}
        & 4.0
        & 8.1
        & 19.5
        & 79.5
        & 262.5
        & 547.8
        & 934.6
        & 1422.6
        %& $O\paren{\frac{L^2\overline{R}^2}{N^2}}$
        \\
\cline{2-10}
        & \eqref{eq:init,func}
        & 4.0
        & 8.1
        & 19.5
        & 79.5
        & 262.5
        & 547.8
        & 934.6
        & 1422.6
        %& $O\paren{\frac{L^2R^2}{N^2}}$
        \\
\cline{2-10}
        & (IFC$'$) %\eqref{eq:init,func}
        & 3.0
        & 7.1
        & 18.5
        & 78.5
        & 261.5
        & 546.8
        & 933.6
        & 1421.6
        %& $O\paren{\frac{L^2R^2}{N^2}}$
	\\
\hline
\end{tabular}
\end{table}

Table~\ref{tab:ratePESTO}
illustrates that the worst-case gradient rates of OGM-G
are equivalent numerically under both~\eqref{eq:init,dist} and~\eqref{eq:init,func}.
This is because the worst-case problem instance of OGM-G in Lemma~\ref{lem:worst} 
associated with the quadratic function
%that is one of the worst-case functions for OGM-G 
under~\eqref{eq:init,func}
also serves as a worst-case of OGM-G under~\eqref{eq:init,dist},
as formally discussed next.

\begin{corollary}
\label{cor}
Let $\x_0,\ldots,\x_N\in\reals^d$ be generated by OGM-G.
Then, for any $N\ge1$,
\begin{align}
\frac{L^2\overline{R}^2}{\ttheta_0^2} 
\le \min_{\substack{f\in\cF, \\ \x_*\in\Xs, \\ ||\x_0-\x_*||\le\overline{R}}}
	||\nabla f(\x_N)||^2
\label{eq:low}
.\end{align}
%with an initial point $\x_0$ satisfying
%$||\x_0 - \x_*|| = \overline{R}$
\begin{proof}
Consider the quadratic function $\phi_2(\x) = \frac{L}{2}||\x||^2$ 
in Lemma~\ref{lem:worst}
associated with the initial point $\x_0=R\nnu$ for any unit-norm vector $\nnu$.
This initial point $\x_0$ satisfies $||\x_0-\x_*||=R$
as well as $\phi_2(\x_0) - \phi_{2,*} = \frac{1}{2}LR^2$,
which implies the inequality~\eqref{eq:low} based on Lemma~\ref{lem:worst}.
\qed
\end{proof}
\end{corollary}

We conjecture that the lower bound~\eqref{eq:low} of OGM-G under~\eqref{eq:init,dist}
is exact, based on numerical evidence in Table~\ref{tab:ratePESTO}.
This is a bit disappointing, 
because it appears that a method that is optimal under 
one initial condition
is far from optimal for another initial condition.
It is also unfortunate that OGM-G
has a poor worst-case rate for decreasing the cost function
under~\eqref{eq:init,dist}.
An open problem is finding a method
that achieves optimal rates invariant to worst-case rate measures
and initial conditions.

In addition,
we study how the worst-case rate %(of OGM-G) 
under~\eqref{eq:init,func}
transfers to that under~\eqref{eq:init,dist}
for given problem instance $(f,\x_0)$.
We particularly focus on two worst-case problem instances of OGM-G in Lemma~\ref{lem:worst},
while similar analysis can be done for the worst-case problem instance of GM
in Lemma~\ref{lem:worst,gm}.
For the worst-case of OGM-G
associated with the Huber function $\phi_1(\x)$,
the constants $R$ and $\overline{R}$ in~\eqref{eq:init,func} and~\eqref{eq:init,dist}
have the following relationship:
\begin{align}
\overline{R} = ||\x_0 - \x_*|| = \frac{\ttheta_0^2+1}{2\ttheta_0}R
\ge \frac{\ttheta_0}{2}R
\ge \frac{N+1}{2\sqrt{2}}R
.\end{align}
%Note that in the worst-case,
%$\overline{R}$ is not really a constant,
%but rather has a value
%that depends on both $R$ and $N$.
We can then show
the following upper bound associated with $\overline{R}$
after $N$ iterations of OGM-G:
\begin{align}
||\nabla \phi_1(\x_N)||^2 = \frac{L^2R^2}{\ttheta_0^2}
        \le \frac{2L^2R^2}{(N+1)^2}
        \le \frac{16L^2\overline{R}^2}{(N+1)^4}
,\end{align}
yielding $O(1/N^4)$, instead of the OGM-G rate $O(1/N^2)$,
expressed by using $\overline{R}$ instead of $R$.
On the other hand, for the worst-case of OGM-G
associated with the quadratic function $\phi_2(\x)$ in Lemma~\ref{lem:worst},
we have the relationship $R=\overline{R}$,
as mentioned in Corollary~\ref{cor}.
These examples illustrate that comparing the worst-case rates
under different initial conditions
is subtle,
and it would be incomplete to treat $R$ and $\overline{R}$
as just arbitrary constants (unrelated to $N$)
in Table~\ref{tab:rate}.

\subsection{\bf Related Work: OGM}
\label{sec:relatedOGM}

This section shows that the proposed OGM-G has a close relationship
with the following OGM~\cite{kim:16:ofo}
(that was numerically first identified in~\cite{drori:14:pof}).

\fbox{
\begin{minipage}[t]{0.85\linewidth}
\vspace{-10pt}
\begin{flalign}
&\quad \text{\bf OGM~\cite{kim:16:ofo}} & \nonumber \\
&\qquad \text{Input: } f\in\cF,\; \x_0=\y_0\in\reals^d,\; N\ge1, \;\htheta_0 = 1. 
        & \nonumber \\
&\qquad \text{For } i = 0,\ldots,N-1, & \nonumber\\
&\qquad \qquad \y_{i+1} = \x_i - \frac{1}{L}\nabla f(\x_i), & \nonumber \\
&\qquad \qquad \htheta_{i+1} = \begin{cases}
                                \frac{1 + \sqrt{1+4\htheta_i^2}}{2}, & i<N-1, \\
                                \frac{1 + \sqrt{1+8\htheta_i^2}}{2}, & i=N-1,
                        \end{cases} & \label{eq:theta} \\
&\qquad \qquad \x_{i+1} = \y_{i+1}
                + \frac{\htheta_i - 1}
                        {\htheta_{i+1}}(\y_{i+1} - \y_i)
                + \frac{\htheta_i}{\htheta_{i+1}}(\y_{i+1} - \x_i). & \nonumber
\end{flalign}
\end{minipage}
} \vspace{5pt}

We can easily notice
the symmetric relationship of the parameters
\begin{align}
\htheta_i = \ttheta_{N-i}, \quad i=0,\ldots,N,
\label{eq:theta,rel}
\end{align}
and the fact that OGM and OGM-G have forms
that differ in the coefficients
of the terms $\y_{i+1} - \y_i$ and $\y_{i+1} - \x_i$.

For $f\in\cF$, $\x_*\in\Xs$ and $||\x_0-\x_*||\le\overline{R}$~\eqref{eq:init,dist},
the final $N$th iterate $\x_N$ of OGM 
has the following exact worst-case cost function bound
\cite[Theorems~2 and 3]{kim:16:ofo}:
\begin{align}
f(\x_N) - f(\x_*) \le \frac{L\overline{R}^2}{2\htheta_N^2}
	\le \frac{L\overline{R}^2}{(N+1)^2}
\label{eq:ogm,rate}
,\end{align}
where this exact upper bound is 
equivalent to the exact worst-case gradient bound~\eqref{eq:ogmg,rate} of OGM-G 
up to a constant $\frac{\overline{R}^2}{2LR^2}$.
This equivalence is similar to the relationship between
the exact worst-case bounds~\eqref{eq:gm_bound_ic2} and~\eqref{eq:gm_cost}
of GM discussed in Remark~\ref{rk:gm}.
The worst-case rate \eqref{eq:ogm,rate} of OGM
is exactly optimal
for large-dimensional smooth convex minimization~\cite{drori:17:tei}.

OGM is equivalent to FSFOM 
with the step coefficients~\cite[Proposition 4]{kim:16:ofo}:
\begin{align}
\hhh_{i+1,k} = \begin{cases}
        \frac{\htheta_i - 1}{\htheta_{i+1}}\hhh_{i,k}, & k=0,\ldots,i-2, \\
        \frac{\htheta_i - 1}{\htheta_{i+1}}(\hhh_{i,i-1} - 1), & k=i-1, \\
        1 + \frac{2\htheta_i - 1}{\htheta_{i+1}}, & k=i.
\end{cases}
\label{eq:ogm,h}
\end{align}
for $i=0,\ldots,N-1$.
The following proposition shows
the symmetric relationship between 
the step coefficients $\{\hhh_{i+1,k}\}$~\eqref{eq:ogm,h} 
and $\{\tth_{i+1,k}\}$~\eqref{eq:ogmg,h}
of OGM and OGM-G, respectively.

\begin{proposition}
%\label{prop:ogmg}
The step coefficients $\{\hhh_{i+1,k}\}$~\eqref{eq:ogm,h}
and $\{\tth_{i+1,k}\}$~\eqref{eq:ogmg,h}
of OGM and OGM-G, respectively, have the following relationship
\begin{align}
\hhh_{i+1,k} = \tth_{N-k,N-i-1}, 
\quad i=0,\ldots,N-1,\;k=0,\ldots,i.
\label{eq:h,rel}
\end{align}
\begin{proof}
We use induction.
Obviously, $\hhh_{1,0} = \tth_{N,N-1}$.
Then,
assuming $\hhh_{i+1,k} = \tth_{N-k,N-i-1}$ for $i=0,\ldots,n-1$,
we have
\begin{align*}
\hhh_{n+1,k} 
%&= \begin{cases}
%        \frac{\htheta_n - 1}{\htheta_{n+1}}\hhh_{n,k}, & k=0,\ldots,n-2, \\
%        \frac{\htheta_n - 1}{\htheta_{n+1}}(\hhh_{n,n-1} - 1), & k=n-1, \\
%        1 + \frac{2\htheta_n - 1}{\htheta_{n+1}}, & k=n,
%	\end{cases}	
&= \begin{cases}
	\frac{\ttheta_{N-n} - 1}{\ttheta_{N-n-1}}\tth_{N-k,N-n}, & k=0,\ldots,n-2, \\
        \frac{\ttheta_{N-n} - 1}{\ttheta_{N-n-1}}(\tth_{N-n+1,N-n} - 1), & k=n-1, \\
        1 + \frac{2\ttheta_{N-n} - 1}{\ttheta_{N-n-1}}, & k=n,
	\end{cases} \\
&= \tth_{N-k,N-n-1}
.\end{align*}
\qed
\end{proof}
\end{proposition}

Building upon the relationships~\eqref{eq:theta,rel} and~\eqref{eq:h,rel}
between OGM and OGM-G,
we numerically study the momentum coefficient values $\beta_i$ and $\gamma_i$
of OGM and OGM-G
in the following form,
to characterize the convergence behaviors of the methods.

\fbox{
\begin{minipage}[t]{0.85\linewidth}
\vspace{-10pt}
\begin{flalign}
&\quad \text{\bf Accelerated First-Order Method} & \nonumber \\
&\qquad \text{Input: } f\in\cF,\; \x_0=\y_0\in\reals^d,\; N\ge1. 
        & \nonumber \\
&\qquad \text{For } i = 0,\ldots,N-1, & \nonumber\\
&\qquad \qquad \y_{i+1} = \x_i - \frac{1}{L}\nabla f(\x_i), & \nonumber \\
&\qquad \qquad \x_{i+1} = \y_{i+1}
                + \beta_i(\y_{i+1} - \y_i)
                + \gamma_i(\y_{i+1} - \x_i). & \nonumber
\end{flalign}
\end{minipage}
} \vspace{5pt}

\begin{figure}[h!]
\begin{center}
\includegraphics[clip,width=0.80\textwidth]{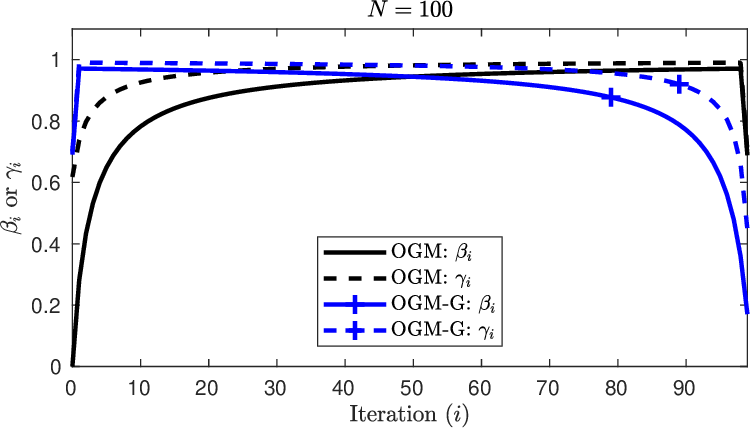}
\end{center}
\caption{Comparison of momentum coefficients ($\beta_i$,$\gamma_i$) of OGM and OGM-G.}
\label{fig:coeff}
\end{figure}

Figure~\ref{fig:coeff} 
compares the momentum coefficients ($\beta_i$, $\gamma_i$)
of OGM and OGM-G for $N=100$.
It is notable that having \emph{increasing} values
of ($\beta_i$, $\gamma_i$) 
as $i$ increases,
except for the last iteration,
yields the optimal (fast) worst-case rate for decreasing the cost function,
whereas having \emph{decreasing} values 
of ($\beta_i$, $\gamma_i$),
except for the first iteration,
yields the fast worst-case rate (that is optimal up to a constant)
for decreasing the gradient.
We leave further theoretical study on such choices of coefficients
as future work.

We next compare OGM and OGM-G
with their other equivalent efficient forms.
Similar to~\cite[Algorithm OGM2]{kim:16:ofo}
one can easily show that the last line of OGM is equivalent to
\begin{align}
\begin{cases}
\z_{i+1} = \y_{i+1} + (\htheta_i - 1)(\y_{i+1} - \y_i)
                + \htheta_i(\y_{i+1} - \x_i), & \\
\x_{i+1} = \paren{1 - \frac{1}{\htheta_{i+1}}}\y_{i+1}
                + \frac{1}{\htheta_{i+1}}\z_{i+1}, &
\end{cases}
\end{align}
while that of OGM-G is equivalent to
\begin{align}
\begin{cases}
\z_{i+1} = \y_{i+1} + (\ttheta_i - 1)(\y_{i+1} - \y_i)
               + \ttheta_i(\y_{i+1} - \x_i), & \\
\x_{i+1} =
        \paren{1 - \frac{2\ttheta_{i+1}-1}{\ttheta_i(2\ttheta_i - 1)}}\y_{i+1}
        + \frac{2\ttheta_{i+1}-1}{\ttheta_i(2\ttheta_i - 1)}\z_{i+1}. &
\end{cases}
\end{align}
This interpretation stems from a variant of FGM~\cite{nesterov:05:smo} 
that involves a convex combination of two points as above.
\cite{kim:16:ofo} already showed that similar interpretation is possible for OGM,
and the expression here also implies that decreasing gradient can be achieved
via some convex combination of two points.
Further analysis is left as future work.

\section{Conclusions}
\label{sec:conc}

This paper developed a first-order method named OGM-G
that has an inexpensive per-iteration computational complexity
and achieves the optimal worst-case bound
for decreasing the gradient of large-dimensional smooth convex functions
up to a constant,
under the initial bounded function condition.
A simple method in~\cite{nesterov:20:pda}, using the OGM-G,
also achieves the optimal worst-case gradient bound
up to a constant,
under the initial bounded distance condition.
The OGM-G was derived by optimizing the step coefficients of first-order methods
in terms of the worst-case gradient bound
using the performance estimation problem (PEP) approach~\cite{drori:14:pof}.
On the way, the exact worst-case gradient bound for a gradient method
was studied.

A practical drawback of OGM-G is that
one must choose the number of iterations $N$ in advance.
Finding a first-order method
that achieves the optimal worst-case gradient bound (up to a constant),
but that does not depend on selecting $N$ in advance,
remains an open problem.
In addition, extending the approaches based on PEP in this paper
to the initial bounded distance condition~\eqref{eq:init,dist}
will be interesting future work;
this PEP approach with a strict relaxation (unlike this paper)
has been studied in~\cite{kim:18:gto}.
Further extensions of this paper to nonconvex problems
and composite problems
are also of interest.

%\section*{Software}
%
%current codes in \verb|/work/ogm_strongly_convex|
%
%Matlab codes
%are available at
%https://gitlab.eecs.umich.edu/michigan-fast-optimization.

\begin{acknowledgements}
Part of this work was carried through while the first author
was affiliated with the University of Michigan.
The first author would like to thank
Ernest K. Ryu
for pointing out related references.
The authors would like to thank associate editor and referees for useful comments,
especially regarding the case where a finite minimizer
does not exist.
The first author was supported in part by the National Research Foundation of Korea (NRF) grant
funded by the Korea government (MSIT) (No. 2019R1A5A1028324),
and the POSCO Science Fellowship of POSCO TJ Park Foundation.
The second author was supported in part by NSF grant IIS 1838179.
\end{acknowledgements}

\appendix

\section*{Appendix: Proof of Eqs.~\eqref{eq:ogmg,h,sum1} and \eqref{eq:ogmg,h,sum2}}
%\label{appx:ogmg,h}

This proof shows the properties~\eqref{eq:ogmg,h,sum1}
and~\eqref{eq:ogmg,h,sum2} 
of the step coefficients $\{\tth_{i,j}\}$~\eqref{eq:ogmg,h}.

We first show~\eqref{eq:ogmg,h,sum1}.
We can easily derive
\begin{align*}
\tth_{i,i-2} 
= \frac{(\ttheta_{i-1}-1)(2\ttheta_i-1)}{\ttheta_{i-2}\ttheta_{i-1}}
= \frac{\ttheta_i^2(2\ttheta_i-1)}{\ttheta_{i-2}\ttheta_{i-1}^2}
\end{align*}
for $i=2,\ldots,N$ 
using~\eqref{eq:ttheta,rule}.
Again using the definition of~\eqref{eq:ogmg,h} and~\eqref{eq:ttheta,rule},
we have
\begingroup
\allowdisplaybreaks
\begin{align*}
\tth_{i,j} &= \frac{\ttheta_{j+1}-1}{\ttheta_j}\tth_{i,j+1}
	%&= \frac{\ttheta_{j+1}-1}{\ttheta_j}
	%	\frac{\ttheta_{j+2}-1}{\ttheta_{j+1}}
	%	\tth_{i,j+2} \\
	%&\qquad\quad\vdots \\
	= \cdots
	= \paren{\prod_{l=j+1}^{i-2}\frac{\ttheta_l-1}{\ttheta_{l-1}}}
                \tth_{i,i-2}
        = \paren{\prod_{l=j+1}^{i-1}\frac{\ttheta_l-1}{\ttheta_{l-1}}}
                \frac{2\ttheta_i-1}{\ttheta_{i-1}} \\
        &= \frac{1}{\ttheta_j}\frac{1}{\ttheta_{j+1}}
                \frac{\ttheta_{j+1}-1}{\ttheta_{j+2}}
                \cdots
                \frac{\ttheta_{i-3}-1}{\ttheta_{i-2}}
                (\ttheta_{i-2}-1)(\ttheta_{i-1}-1)
                \frac{2\ttheta_i-1}{\ttheta_{i-1}} \\
        &= \frac{1}{\ttheta_j}\frac{1}{\ttheta_{j+1}}
                \frac{\ttheta_{j+2}}{\ttheta_{j+1}}
                \cdots
                \frac{\ttheta_{i-2}}{\ttheta_{i-3}}
                (\ttheta_{i-2}-1)(\ttheta_{i-1}-1)
                \frac{2\ttheta_i-1}{\ttheta_{i-1}} \\
        &= \frac{\ttheta_{i-2}(\ttheta_{i-2}-1)(\ttheta_{i-1}-1)(2\ttheta_i-1)}
                {\ttheta_j\ttheta_{j+1}^2\ttheta_{i-1}}
        = \frac{\ttheta_i^2(2\ttheta_i-1)}{\ttheta_j\ttheta_{j+1}^2}, 
\end{align*}
\endgroup
for $i=2,\ldots,N,\;j=0,\ldots,i-3$,
which concludes the proof of~\eqref{eq:ogmg,h,sum1}.

We next prove the first two lines of~\eqref{eq:ogmg,h,sum2} using the induction.
For $N=1$, we have $\ttheta_1 = 1$ and 
\begin{align*}
\tth_{1,0} = 1 + \frac{2\ttheta_1-1}{\ttheta_0}
	= 1 + \frac{\ttheta_1^2}{\ttheta_0}
	= 1 + \frac{\frac{1}{2}(\ttheta_0^2 - \ttheta_0)}{\ttheta_0}
	= \frac{1}{2}(\ttheta_0+1)
,\end{align*}
where the third equality uses~\eqref{eq:ttheta,rule}.
For $N>1$, we have 
\begin{align*}
\tth_{N,N-1} = 1 + \frac{2\ttheta_N-1}{\ttheta_{N-1}} 
	= 1 + \frac{\ttheta_N^2}{\ttheta_{N-1}} 
	= 1 + \frac{\ttheta_{N-1}^2 - \ttheta_{N-1}}{\ttheta_{N-1}}
	= \ttheta_{N-1}
,\end{align*}
where the third equality uses~\eqref{eq:ttheta,rule}.
Assuming $\sum_{l=j+1}^N\tth_{l,j} = \ttheta_j$
for $j=n,\ldots,N-1$ and $n\ge1$, we get
\begingroup
\allowdisplaybreaks
\begin{align*}
\sum_{l=n}^N\tth_{l,n-1}
=\;& 1 + \frac{2\ttheta_n-1}{\ttheta_{n-1}} 
	+ \frac{\ttheta_n-1}{\ttheta_{n-1}}(\tth_{n+1,n}-1)
	+ \frac{\ttheta_n-1}{\ttheta_{n-1}}\sum_{l=n+2}^N\tth_{l,n} \\
=\;& 1 + \frac{\ttheta_n}{\ttheta_{n-1}}
	+ \frac{\ttheta_n-1}{\ttheta_{n-1}}\sum_{l=n+1}^N\tth_{l,n}
= \frac{\ttheta_{n-1} + \ttheta_n + (\ttheta_n-1)\ttheta_n}{\ttheta_{n-1}}
= \frac{\ttheta_{n-1} + \ttheta_n^2}{\ttheta_{n-1}} \\
=\;& \begin{cases}
\frac{1}{2}(\ttheta_0 + 1), & n = 0, \\
\ttheta_n, & n=1,\ldots,N-1,
\end{cases}
\end{align*}
\endgroup
where the last equality uses~\eqref{eq:ttheta,rule},
which concludes the proof of the first two lines of~\eqref{eq:ogmg,h,sum2}.

We finally prove the last line of~\eqref{eq:ogmg,h,sum2} using the induction.
For $i\ge1$, we have
\begin{align*}
\sum_{l=i+1}^N\tth_{l,i-1}
	&= \sum_{l=i}^N\tth_{l,i-1} - \tth_{i,i-1}
	= \ttheta_{i-1} - \paren{1+\frac{2\ttheta_i-1}{\ttheta_{i-1}}}
	= \frac{(\ttheta_i-1)^2}{\ttheta_{i-1}}
	= \frac{\ttheta_{i+1}^4}{\ttheta_{i-1}\ttheta_i^2} 
,\end{align*}
where the third and fourth equalities use~\eqref{eq:ttheta,rule}.
Then, assuming $\sum_{l=i+1}^N\tth_{l,j}=\frac{\ttheta_i^4}{\ttheta_j\ttheta_{j+1}^2}$ 
for $i=n,\ldots,N-1$, $j=0,\ldots,i-1$ with $n\ge1$,
we get:
\begin{align*}
\sum_{l=n}^N\tth_{l,j} 
	&= \sum_{l=n+1}^N\tth_{l,j} + \tth_{n,j}
	= \frac{\ttheta_{n+1}^4}{\ttheta_j\ttheta_{j+1}^2} 
		+ \frac{\ttheta_n^2(2\ttheta_n-1)}{\ttheta_j\ttheta_{j+1}^2}
	= \frac{\ttheta_n^2(\ttheta_n-1)^2 + \ttheta_n^2(2\ttheta_n-1)}{\ttheta_j\ttheta_{j+1}^2}
	= \frac{\ttheta_n^4}{\ttheta_j\ttheta_{j+1}^2}
,\end{align*}
where the second and third equalities use~\eqref{eq:ogmg,h,sum1},
which concludes the proof.
\qed

%\ifjota
	\bibliographystyle{spmpsci_unsrt}
%\else
	%\bibliographystyle{spmpsci}
%\fi
\bibliography{master,mastersub}
\end{document}